\documentclass[final]{siamltex1213}
\usepackage{latexsym,amssymb,amsmath,amsfonts,float}
\usepackage{algorithmicx}
\usepackage{algpseudocode}
\usepackage{tikz}

\usepackage{bm}
\usepackage{graphicx}
\usepackage{stmaryrd}

\usepackage{enumerate}

\def\CC{{\mathbb C}}

\def\NN{{\mathbb N}}

\def\RR{{\mathbb R}}

\def\cal{\mathcal}

\def\cH{{\cal H}}

\def\cO{{\cal O}}

\def\g{{\mathbf{g}}}
\def\w{{\mathbf{w}}}

\newcommand{\C}{\ensuremath{\mathbb{C}}}

\newtheorem{algorithm}{Algorithm}

\def\C{\mathbf{C}}

\renewcommand{\v}{\mathbf{v}}
\newcommand{\I}{\mathbf{I}}
\renewcommand{\H}{\mathbf{H}}
\renewcommand{\u}{\mathbf{u}}
\newcommand{\f}{\mathbf{f}}

\newcommand{\x}{\mathbf{x}}

\renewcommand{\k}{\mathbf{k}}

\newcommand{\beq}{\begin{equation}}
\newcommand{\eeq}{\end{equation}}

\newcommand{\bit}{\begin{itemize}}
\newcommand{\eit}{\end{itemize}}
\newcommand{\ben}{\begin{enumerate}}
\newcommand{\een}{\end{enumerate}}
\title{Fast alternating bi-directional preconditioner for the 2D high-frequency Lippmann-Schwinger equation}

\author{Leonardo Zepeda-N\'u\~nez\footnotemark[2] \, and Hongkai Zhao\footnotemark[2]}

\date{August 2015}

\begin{document}

\maketitle

\renewcommand{\thefootnote}{\fnsymbol{footnote}}

\footnotetext[2]{Department of Mathematics, University of California, Irvine, CA 92612, USA.}

\renewcommand{\thefootnote}{\arabic{footnote}}

%\slugger{sisc}{xxxx}{xx}{x}{x--x}%slugger should be set to mms, siap, sicomp, sicon, sidma, sima, simax, sinum, siopt, sisc, or sirev

\begin{abstract}
This paper presents a fast iterative solver for Lippmann-Schwinger equation for high-frequency waves scattered by a smooth medium with a compactly supported inhomogeneity. The solver is based on the sparsifying preconditioner \cite{Ying:Sparsifying_Preconditioner_for_the_Lippmann--Schwinger_Equation} and a domain decomposition approach similar to the method of polarized traces \cite{ZepedaDemanet:the_method_of_polarized_traces}. The iterative solver has two levels, the outer level in which a sparsifying preconditioner for the Lippmann-Schwinger equation is constructed, and the inner level, in which the resulting sparsified system is solved fast using an iterative solver preconditioned with a bi-directional matrix-free variant of the method of polarized traces. The complexity of the construction and application of the preconditioner is $\cO(N)$ and $\cO(N\log{N})$ respectively, where $N$ is the number of degrees of freedom. Numerical experiments in 2D indicate that the number of iterations in both levels depends weakly on the frequency resulting in a method with an overall $\cO(N\log{N})$ complexity.
\end{abstract}

\section{Introduction}
In this paper we present a fast solver for the high-frequency Lippmann-Schwinger equation in 2D. Let $\omega$ be the frequency of a time-harmonic wave propagating in a medium with wave-speed denoted by $c(\x)$, where $\x = (x,z)$. Let $\Omega$ be a rectangular domain and assume that the refraction index is constant up to a compactly supported perturbation. The refraction index is denoted by $1/c^2(\x) = 1 + m(\x)$, where $\mbox{supp}(m) = D \subset \Omega \subset \RR^2$.
In this case, $u(\x)$, the wave scattered by the perturbation $m$ from an incident wave $u_{I}(\x)$, is the solution to the Helmholtz equation
\begin{equation} \label{eq:Helmholtz}
	\left ( \Delta + \omega^2( 1 + m(\x) ) \right ) (u(\x) + u_{I}(\x)) = 0, \qquad \x \in \RR^d;
\end{equation}
where $u(\x)$ satisfies the Sommerfeld radiation condition,
\begin{equation}
	\lim_{ r \rightarrow \infty} r^{1/2} \left (  \frac{\partial u}{\partial r} - i\omega u\right) = 0,
\end{equation}
and the incoming field $u_{I}(\x)$ satisfies the homogeneous Helmholtz equation
\begin{equation} \label{eq:Helmholtz_incoming}
	\left ( \Delta + \omega^2 \right )u_{I}(\x) = 0,
\end{equation}
in some neighborhood of $\Omega$ containing $D$.
% %%% We keep this for higher dimensional problems
% For the free space Helmholtz operator, $-\Delta - \omega^2$ with Sommerfeld radiation conditions, the Green's function is given by
% \begin{equation} \label{eq:Greens_function}
% 	G(\x) = \left \{ \begin{array}{ll}
% 					\frac{1}{2i\omega} e^{i\omega |\x|}  & , d=1 \\
% 					\frac{i}{4} H^1_0(\omega |\x|)       & , d=2 \\
% 					\frac{1}{4\pi |\x|} e^{i\omega |x|}  & , d=3 \\
% 					 \end{array}
% 			\right .
% \end{equation}

From eqs. \eqref{eq:Helmholtz} and \eqref{eq:Helmholtz_incoming} it follows that the unknown field $u$ satisfies
\begin{equation} \label{eq:Helmholtz_2}
	\left ( \Delta + \omega^2( 1 + m(\x) ) \right ) u(\x)  = -\omega^2 m(\x) u_{I}(\x), \qquad \x \in \RR^d.
\end{equation}
Convolving \eqref{eq:Helmholtz_2} with
\begin{equation} \label{eq:Greens_function}
	G(\x) = -\frac{i}{4} H^1_0(\omega |\x|),
\end{equation}
the Green's function for the 2D free-space Helmholtz operator, results in the Lippmann-Schwinger equation:
\begin{equation} \label{eq:Lippmann-Schwinger}
	u + \omega^2 G*(m u) = -\omega^2 G*(m u_{I}).
\end{equation}

We chose \eqref{eq:Lippmann-Schwinger} among different formulations of the Lippmann-Schwinger equation (for other formulations see \cite{ambikasaran_greengard:Fast_adaptive_high_order_accurate_discretization_of_the_Lippmann-Schwinger_equation_in_two_dimension}) by the following reasons :

\begin{itemize}
\item the unknown scattered field, $u(\x)$, satisfies the Sommerfeld radiation conditions automatically;
\item the convolution kernel is translation invariant which allows for a rapid application via fast Fourier transforms (FFT) \cite{Cooley:An_Algorithm_for_the_Machine_Calculation_of_Complex_Fourier_Series}.
\end{itemize}

%Even though $u$ can be computed by discretizing \eqref{eq:Helmholtz_2} using finite differences or finite elements; and then solving the resulting system using state-of-the art preconditioners \cite{Chen_Xiang:a_source_transfer_ddm_for_helmholtz_equations_in_unbounded_domain,EngquistYing:Sweeping_H,EngquistYing:Sweeping_PML,CStolk_rapidily_converging_domain_decomposition,GeuzaineVion:double_sweep,ZepedaDemanet:the_method_of_polarized_traces}, there are several advantages on computing $u$ by solving \eqref{eq:Lippmann-Schwinger} :
Even though $u$ can be computed by discretizing \eqref{eq:Helmholtz_2} using finite differences, finite elements or composite spectral discretizations \cite{Gillman_Barnett_Martinsson:A_spectrally_accurate_solution_technique_for_frequency_domain_scattering_problems_with_variable_media,Gillman:A_Direct_Solver_with_ON_Complexity_for_Variable_Coefficient_Elliptic_PDEs_Discretized_via_a_High_Order_Composite_Spectral_Collocation_Method}; and then solving the resulting linear system using state-of-the art preconditioners \cite{Chen_Xiang:a_source_transfer_ddm_for_helmholtz_equations_in_unbounded_domain,EngquistYing:Sweeping_H,EngquistYing:Sweeping_PML,CStolk_rapidily_converging_domain_decomposition,GeuzaineVion:double_sweep,ZepedaDemanet:the_method_of_polarized_traces}, there are several advantages on computing $u$ by solving \eqref{eq:Lippmann-Schwinger} :
\begin{enumerate}
\item \eqref{eq:Lippmann-Schwinger} arises from a second kind Fredholm formulation, implying that the conditioning number of the discretized system does not deteriorate as the discretization is refined;

\item \eqref{eq:Lippmann-Schwinger} is defined on a bounded set $\Omega$ and the solution satisfies the Sommerfeld radiation condition automatically, whereas \eqref{eq:Helmholtz_2} is defined in the whole $\RR^2$; hence one needs to truncate the computational domain, and to impose absorbing boundary conditions (ABC) at the boundary of the truncated computational domain, using perfectly matched layers (PML) \cite{Berenger:PML,Bermudez:An_optimal_perfectly_matched_layer_with_unbounded_absorbing_function_for_time-harmonic_acoustic_scattering_problems,Johnson:PML} or other techniques \cite{Engquist:Absorbing_boundary_conditions_for_the_numerical_simulation_of_waves} (see \cite{Bermudez:Perfectly_Matched_Layers_for_Time-Harmonic_Second_Order_Elliptic_Problems} for a review);

% \item  most discretizations of \eqref{eq:Helmholtz_2} via finite differences or finite elements suffer from the pollution effect (\cite{Sauter_Babuska:Is_the_Pollution_Effect_of_the_FEM_Avoidable_for_the_Helmholtz_Equation_Considering_High_Wave_Numbers,Bayliss:On_accuracy_conditions_for_the_numerical_computation_of_waves,Babuska:A_Generalized_Finite_Element_Method_for_solving_the_Helmholtz_equation_in_two_dimensions_with_minimal_pollution} (see  \cite{Stolk,Turkel:Compact_2D_and_3D_sixth_order_schemes_for_the_Helmholtz_equation_with_variable_wave_number} for recent references to methods designed to alleviate this issue); whereas the introduction of the analytical Green's function in \eqref{eq:Lippmann-Schwinger} automatically alleviate the pollution.
\item  most discretizations of \eqref{eq:Helmholtz_2} via finite differences or finite elements suffer from dispersion error and the pollution effect \cite{Babuska:A_Generalized_Finite_Element_Method_for_solving_the_Helmholtz_equation_in_two_dimensions_with_minimal_pollution,Sauter_Babuska:Is_the_Pollution_Effect_of_the_FEM_Avoidable_for_the_Helmholtz_Equation_Considering_High_Wave_Numbers,Bayliss:On_accuracy_conditions_for_the_numerical_computation_of_waves}  (see  \cite{Stolk,Turkel:Compact_2D_and_3D_sixth_order_schemes_for_the_Helmholtz_equation_with_variable_wave_number} for recent references to methods designed to alleviate this issue); whereas the introduction of the analytical Green's function in \eqref{eq:Lippmann-Schwinger} and the use of an accurate quadrature rule automatically alleviates these issues.

\end{enumerate}

However, solving \eqref{eq:Lippmann-Schwinger} raises other numerical issues. The discretization of \eqref{eq:Lippmann-Schwinger} produces large dense linear systems, which becomes impractical, both memory- and complexity-wise, for typical direct solvers. Although a new generation of  direct methods based on hierarchical compression, such as \cite{ambikasaran_greengard:Fast_adaptive_high_order_accurate_discretization_of_the_Lippmann-Schwinger_equation_in_two_dimension,Cheng:An_adaptive_fast_solver_for_the_modified_Helmholtz_equation_in_two_dimensions,Corona:An_direct_solver_for_integral_equations_on_the_plane,Ho:Hierarchical_Interpolative_Factorization_for_Elliptic_Operators_Integral_Equations}, have been able to considerably reduce the complexity of solving integral equations coming from elliptic partial differential equations (PDE),  when applying these methods to integral equations coming from wave equations such as \eqref{eq:Lippmann-Schwinger}, they manage to exhibit linear or quasi-linear complexity for low frequency problems, but they fall back to $\cO(N^{3/2})$ in high-frequency regimes. This phenomenon can be clearly explained by recent studies of the approximate separability of the Green's functions.  It was shown in \cite{Bebendorf:Existence_of_Hmatrix_approximants_to_the_inverse_FE_matrix_of_elliptic_operators_with_Linftycoefficients} that
the Green's functions for elliptic operators are highly separable, implying that the off-diagonal blocks of the discrete Green's function are extremely low-rank. In contrast, it was shown in \cite{Engquist_Zhao:approximate_separability_of_green_function_for_high_frequency_Helmholtz_equations} that the Green's functions for the Helmholtz operator in high-frequency regime are highly non-separable in general. This last results  provides an asymptotic frequency-dependent lower bound on the ranks of the off-diagonal blocks of the discrete Green's function, which hinders the compressibility of the blocks, thus resulting in algorithms with higher complexity.

In addition, the Lippmann-Schwinger equation becomes numerically-ill conditioned as the frequency and contrast increases, due to multiple scattering, hindering the convergence rate of iterative methods. This difficulty has led to the development of efficient preconditioners, such as \cite{Andersson:A_Fast_Bandlimited_Solver_for_Scattering_Problems_in_Inhomogeneous_Media,Bruno:Wave_scattering_by_inhomogeneous_media_efficient_algorithms_and_applications,Bruno:An_efficient_preconditioned_high_order_solver_for_scattering_by_two-dimensional_inhomogeneous_media,Sifuntes:Preconditioned_Iterative_Methods_for_Inhomogeneous_Acoustic_Scattering_Applications,Vainikko:Fast_Solvers_of_the_Lippmann-Schwinger_Equation}; however, they tend to require significantly more iterations to converge as the frequency increases. Some of them are based on semi-analytical techniques, which require fairly strict assumptions on the perturbation (for example \cite{Bruno:An_efficient_preconditioned_high_order_solver_for_scattering_by_two-dimensional_inhomogeneous_media}); they tend to work well under the assumptions, but they naturally deteriorate if not.

There is an extensive literature on volume integral methods, including fast methods \cite{Bebendorf:Hierarchical_LU_Decomposition-based_Preconditioners_for_BEM,Ho_Greengard:A_Fast_Direct_Solver_for_Structured_Linear_Systems_by_Recursive_Skeletonization,Martinsson:A_fast_direct_solver_for_boundary_integral_equations_in_two_dimensions} and, in particular, methods solving the Lippmann-Schwinger equation \cite{Amar:Numerical_solution_of_the_Lippmann-Schwinger_equations_in_photoemission:_application_to_xenon,Andersson:A_Fast_Bandlimited_Solver_for_Scattering_Problems_in_Inhomogeneous_Media,Bruno:Wave_scattering_by_inhomogeneous_media_efficient_algorithms_and_applications,Bruno:An_efficient_preconditioned_high_order_solver_for_scattering_by_two-dimensional_inhomogeneous_media,Cheng:An_adaptive_fast_solver_for_the_modified_Helmholtz_equation_in_two_dimensions,Corona:An_direct_solver_for_integral_equations_on_the_plane,Duan-Rohklin:High-order_quadratures_for_the_solution_of_scattering_problems_in_two_dimensions,Lanzara:Numerical_Solution_of_the_Lippmann_Schwinger_Equation_by_Approximate_Approximations,Vainikko:Fast_Solvers_of_the_Lippmann-Schwinger_Equation,Ying:Sparsifying_Preconditioner_for_the_Lippmann--Schwinger_Equation}, and we do not review it here, but to the authors' knowledge there is no method that can handle general smooth perturbations in quasi-linear time in the high-frequency regime.

In the present paper we introduce a new preconditioner for the 2D high-frequency Lippmann-Schwinger equation, with a total asymptotic cost $\cO(N \log{N})$, which compares favorably with respect to most of the methods (both iterative and direct) in the literature. The scaling holds for general smooth perturbations without large resonant cavities. The solver is based on the sparsifying preconditioner and a fast iterative solver in a domain decomposition setting.  In particular, the sparsifying preconditioner is used to approximate the Lippmann-Schwinger equation with a sparse system, which is solved fast using a bi-directional matrix-free variant of the method of polarized traces \cite{Zepeda_Demanet:Nested_domain_decomposition_with_polarized_traces_for_the_2D_Helmholtz_equation}.

The method has two stages, the off-line, or setup, stage that is only performed once for each linear system to solve with an asymptotic cost of $\cO(N)$, and the on-line, or solve, stage that is performed for each right-hand-side, with an asymptotic cost $\cO(N\log{N})$. Furthermore, the on-line stage has two levels: the inner level in which the approximate sparse system is solved iteratively using a bi-directional preconditioner, and the outer level, in which the resulting solution from the inner level is used to precondition the Lippmann-Schwinger equation.

Under the assumptions that the perturbation is smooth and it does not represent a resonant cavity the number of outer iterations is essentially independent of the frequency \cite{Ying:Sparsifying_Preconditioner_for_the_Lippmann--Schwinger_Equation}, and numerical experiments presented in this paper, indicate that the number of inner iterations depends, at worst, logarithmically with the frequency, resulting in the aforementioned complexity.

\subsection{Results}

In this paper we propose a new two-level preconditioner for the 2D high-frequency Lippmann-Schwinger equation that results in a method with $\cO(N \log{N})$ complexity. We follow the approach in \cite{Ying:Sparsifying_Preconditioner_for_the_Lippmann--Schwinger_Equation} to construct a sparsifying preconditioner, and we develop a fast iterative solver to solve the sparsified system, in order to achieve an overall optimal complexity. The main improvement with respect to the sparsifying preconditioner is how the resulting sparse linear system is solved. In \cite{Ying:Sparsifying_Preconditioner_for_the_Lippmann--Schwinger_Equation} the system is solved via multi-frontal methods \cite{Duff_Reid:The_Multifrontal_Solution_of_Indefinite_Sparse_Symmetric_Linear,GeorgeNested_dissection} incurring an asymptotic cost $\cO(N^{3/2})$ for the setup and $\cO(N\log{N})$ for the application. Exchanging the multifrontal solver by an iterative one, preconditioned using a bi-directional variant of the method of polarized traces, results in a drastic reduction of the asymptotic cost.

The inner-level preconditioner is a matrix-free variant of the method of polarized traces \cite{Zepeda_Demanet:Nested_domain_decomposition_with_polarized_traces_for_the_2D_Helmholtz_equation}, which has been adapted to solve the sparsified Lippmann-Schwinger equation efficiently. As it will be explained in the sequel, preconditioning the sparsified system using the original method of polarized traces results in an algorithm that is sensitive to the sign of the perturbation and to the direction of the incident waves. Then a new bi-directional variant is introduced to eliminate these adverse effects.

The method of polarized traces is a domain decomposition method \cite{ZepedaDemanet:the_method_of_polarized_traces}, which relies on a layered domain decomposition coupled with an equivalent surface integral equation (SIE) that is easy to precondition. The matrix-free variant introduced in \cite{Zepeda_Demanet:Nested_domain_decomposition_with_polarized_traces_for_the_2D_Helmholtz_equation} uses local solves on extended subdomains in order to apply the SIE and the preconditioner. The number of iterations to convergence depends on the underlying medium and the quality of the absorbing boundary conditions (ABC) imposed at the interfaces between subdomains. For the Helmholtz problem, there exist a myriad of different formulations of ABC \cite{Berenger:PML,Bermudez:Perfectly_Matched_Layers_for_Time-Harmonic_Second_Order_Elliptic_Problems,Engquist:Absorbing_boundary_conditions_for_the_numerical_simulation_of_waves,Johnson:PML}; unfortunately, for the sparse system resulting from sparsifying the Lippmann-Schwinger equation such ABC are not readily available.
We give details on building such ABC for the sparsified Lippmann-Schwinger system  in
Section \ref{section:ABC}.
% is consecrated to build such ABC for the sparsified Lippmann-Schwinger system.

Alas, the formulation of the ABC in this paper is not perfect, which is evidenced by numerical tests showing that the number of iteration needed to convergence using the matrix-free version of the method of polarized traces grows sub-optimally as $\cO(\sqrt{\omega})$. Optimality can be restored by introducing an alternating bi-directional sweep, and by preconditioning the volume problem directly, which reduces the number of iterations to $\cO(\log{\omega})$.
To show the effectiveness of our proposed method, we also include an example arising from plasma physics, in particular, an electromagnetic wave impinging on a confined plasma in a fusion reactor.

In the present work we mainly focus on explaining the key ideas and verification of the method by 2D tests. The algorithm presented in this paper can be easily extended to 3D, with an offline complexity of $\cO(N^{4/3})$ and an on-line complexity of $\cO(N \log^2{N})$. The complexities can be further reduced using a nested approach such as the one presented in \cite{Liu_Ying:Recursive_sweeping_preconditioner_for_the_3d_helmholtz_equation,Zepeda_Demanet:Nested_domain_decomposition_with_polarized_traces_for_the_2D_Helmholtz_equation} resulting in $\cO(N)$ and $\cO(N \log{N})$ complexity respectively.

\subsection{Outline of the paper}
The present paper is organized as follows: we present the discretization of the Lippmann-Schwinger equation and the provide a brief review of the sparsifying preconditioner in Section \ref{Section:sparsifying_preconditioner}.
In Section 3, we review the method of polarized traces and we discuss the modifications needed to adapt it to the sparsified Lippmann-Schwinger equation, in particular we explain the formulation for the ABC in Section \ref{section:ABC} and the bi-directional preconditioner in Section \ref{section:bidirectional_algorithm}. In Section \ref{section:preconditioner} we present the two-level preconditioner and we discuss its complexity.
Finally, in Section \ref{section:numerical_experiments}, we present various numerical experiments to support the complexity claims.

\section{Discretization and sparsifying preconditioner} \label{Section:sparsifying_preconditioner}

In this section we discuss the discretization used and we provide a brief review of the sparsifying preconditioner \cite{Ying:Sparsifying_Preconditioner_for_the_Lippmann--Schwinger_Equation}. We refer the interested reader to Section 2 of \cite{Ying:Sparsifying_Preconditioner_for_the_Lippmann--Schwinger_Equation} for further details.

We discretize the domain $\Omega$ with an uniform mesh\footnote{If the perturbation has sharp transitions the sparsifying preconditioner can be modified to handle adaptive quadratures without asymptotic penalty. However, the convolution product will have to be handled differently \cite{Beylkin:Fast_convolution_with_the_free_space_Helmholtz_Greens_function} and the construction of the preconditioner will require randomized techniques in order to keep the same asymptotic complexity.} of a step size $h$, such that we have 10 points per wavelength in the background wavespeed, $h\omega = \cO(1)$. We denote the collocations points by $\x_{i,j} = (x_i, z_j) = (hi, hj)$ for $1 \leq i \leq n_x$ and  $1 \leq j \leq n_z$, where $n_x$ and $n_z$ are the number of grid-points in each direction. In addition, we define a vector index $\k = (i,j)$, notation that will be useful to explain the sparsifying preconditioner.
We use a simple modification of the trapezoidal rule \cite{Duan-Rohklin:High-order_quadratures_for_the_solution_of_scattering_problems_in_two_dimensions} to discretize the integral within the convolution operation in \eqref{eq:Lippmann-Schwinger}. In particular, we use a quadrature that only modifies the weight in the diagonal, obtaining a quadrature rule with order $\cO(h^4 \log{h})$.

We denote the unknown scattered field in $\Omega$ as $\u_{\k} = \u_{i,j} = u(\x_{i,j})$, and we write the resulting discretized linear system as
\begin{equation} \label{eq:discrete_Lipmann-Schwinger}
	\H \u = \f,
\end{equation}
in which $\H$ can be written as
\begin{equation}
	\H = \I + \omega^2 \mathbf{G} \cdot \mathbf{D},
\end{equation}
where $\I$ is the identity matrix, and $\mathbf{D}$ is a diagonal matrix that encodes the multiplication of $\u$ times $m(\x_{i,j})$. Finally, $\f$ is given by
\begin{equation}
	\f = -\omega^2 \mathbf{G} \cdot \mathbf{D} \u_{I},
\end{equation}
with $ (\u_{I})_{i,j} = u_{I}(\x_{i,j})$.

Given that the convolution kernel is translation invariant and that the mesh is uniform we can easily embed $\mathbf{G}$ in a Toeplitz matrix and use the FFT to apply it in $\cO(N\log{N})$ operations (or we can use the method in \cite{Beylkin:Fast_convolution_with_the_free_space_Helmholtz_Greens_function} if the mesh is non-uniform). %Given that $\I$ and $\mathbf{D}$ are sparse, it follows that we can apply $\H$ to any vector in $\cO(N\log{N})$ operations.

The sparsifying preconditioner is based on designing a sparse matrix $\mathbf{A}$ such that the product $\mathbf{A} \H$ is approximately sparse in the sense that a large amount of it entries are extremely small. The main idea is to multiply  $\mathbf{A}$ by $\H$, and threshold all the small entries of the product to zero, obtaining a sparse matrix $ \mathbf{C} $ that approximates $\mathbf{A} \H $. As it will be explained in the sequel, the sparsity pattern of $\mathbf{A}$ is the same as a second-order compact finite difference discretization of the Laplacian, and its entries are specially constructed to annihilate the free space Green's function away from the source.  In return, the matrix $\mathbf{C}$, which shares the same sparsity pattern as $\mathbf{A}$, is similar to the matrix arising from a finite differences discretization of the Helmholtz equation \eqref{eq:Helmholtz}.

Since
%\begin{equation}
$	\mathbf{A} \H = \mathbf{A} + \omega^2 \mathbf{A} \mathbf{G} \mathbf{D}$,
%\end{equation}
 it follows that to make $\mathbf{A} \mathbf{H}$ approximately sparse, we only need to find $\mathbf{A}$ such that the product  $\mathbf{A} \mathbf{G}$ is approximately sparse. The last requirement is equivalent to enforce that the entries of the matrix $\mathbf{A}$  approximately annihilate the off-diagonal blocks of $\mathbf{G}$. This requirement can be recast as a minimization problem, which is solved via a singular value decomposition (SVD).

Define the set of immediate neighbors (stencils) of $\k$ as
\begin{equation}
	\mu(\k) = \{ \mathbf{p} :  \| \k - \mathbf{p} \| _{\ell_{\infty}(\NN^2)} \leq 1, \x_{\mathbf{p}} \in \Omega \}.
\end{equation}

We seek $\mathbf{A}$ such that for each $\k$:
\begin{itemize}
	\item only $\mathbf{A}(\k, \mu(\k))$ are non-zeros (sparsity of $\mathbf{A}$);
	\item $\mathbf{A}(\k, :)\mathbf{G}(:, \mu(\k)^{c}) \approx 0$ (approximate sparsity of the $\mathbf{A}\mathbf{G}$).
\end{itemize}

Given that $\mathbf{G}$ is translation invariant we impose $\mathbf{A}$ to be translation invariant. In such case, it is sufficient to consider only the case when $\k = 0$. Then the problem can be recast as
\begin{equation}
	\min_{\alpha: \|\alpha \| =1} \|\alpha \cdot \mathbf{G}(\mu(0), \mu(0)^c )  \|,
\end{equation}
which can be efficiently solved via a singular value decomposition of $\mathbf{G}(\mu(0), \mu(0)^c ) = U \Sigma V^*$ such that
\begin{equation}
	\alpha = U(:, \mu(0))^*.
\end{equation}
Finally, $\mathbf{A}$ is built such that $\mathbf{A}(\k, \mu(\k)) = \alpha$, for $\k$ in the interior of $\Omega$. For $\k$ on $\partial \Omega$ a special treatment is needed, which depends on the geometry of $\Omega$. For the sake of simplicity we omit such treatment but we refer the interested reader to Section 4 of \cite{Ying:Sparsifying_Preconditioner_for_the_Lippmann--Schwinger_Equation}.

By multiplying \eqref{eq:discrete_Lipmann-Schwinger} by $\mathbf{A}$ and sparsifying the product, we obtain the approximated {\em sparse} system
\begin{equation}
	\mathbf{C} \u \approx \mathbf{A} \f,
\end{equation}
for which the preconditioner is defined as
\begin{equation} \label{eq:preconditioner}
	\u_{approx} = \mathbf{C}^{-1} \g, \qquad \mbox{where } \g =   \mathbf{A} \f .
\end{equation}

It has been shown in \cite{Ying:Sparsifying_Preconditioner_for_the_Lippmann--Schwinger_Equation} that when the preconditioner mentioned above is used within GMRES iterations, the number of iterations to convergence is almost independent of the frequency.

One alternative to this preconditioner would be to discretize \eqref{eq:Helmholtz_2} via finite differences or finite elements, and then use this approximate solution to precondition \eqref{eq:discrete_Lipmann-Schwinger}. Unfortunately, such preconditioners tend to behave poorly, due to the pollution error that is normally introduced with such discretizations. The pollution error creates serious phase errors making the preconditioner incompatible with the linear system being preconditioned (see \cite{Stolk} for a similar example for multigrid). The optimization procedure used in the construction of the sparsifying preconditioner takes into account the exact free space Green's function and hence minimizes the dispersion error. The size of the stencils used in $\mathbf{A}$ and $\mathbf{C}$ should strike a balance between sparsity and dispersion error.
Since the matrix $\mathbf{A}$ mimics the discretization of the Helmholtz differential operator, it is expected that the ideas for fast solver for the Helmholtz equation can be applied to the resulting sparse system. This is what we propose in the next Section.

\section{Method of polarized traces}
\label{sec:polarized-traces}
As seen from Section  \ref{Section:sparsifying_preconditioner}, the sparsifying preconditioner for the Lippmann-Schwinger equation relies on solving the sparsified system
\begin{equation} \label{eq:global_sparse_system}
	 \mathbf{C}\v = \g,
\end{equation}
which, in the present paper, is solved fast using a bi-directional variant of the method of polarized traces.
In this section, we first briefly introduce the matrix-free version of the method of polarized traces \cite{Zepeda_Demanet:Nested_domain_decomposition_with_polarized_traces_for_the_2D_Helmholtz_equation} in one direction, and then we discuss the modifications needed to adapt it to the sparsified Lippmann-Schwinger equation. We provide a simple description of the three key elements of the method of polarized traces for the Helmholtz equation in the following order:
\begin{enumerate}
	\item  the layered domain decomposition with the extended subdomains;
	\item  how to define the transmission conditions between subdomains using a discrete Green's representation formula issued from an algebraic decomposition;
	\item  how to implement an absorbing boundary conditions for the subdomains to obtain a fast converging iterative method.
\end{enumerate}

We then provide the full algorithm for the one-directional preconditioner; we provide its physical intuition and discuss its short-comings, which are rectified by the bi-directional preconditioner that is introduced right after.

\subsection{Domain decomposition}

\begin{figure}[h]
\begin{center}
\vspace{-0.25cm}
\includegraphics[trim= 5mm 30mm 0mm 20mm, angle =0,clip, width=8cm]{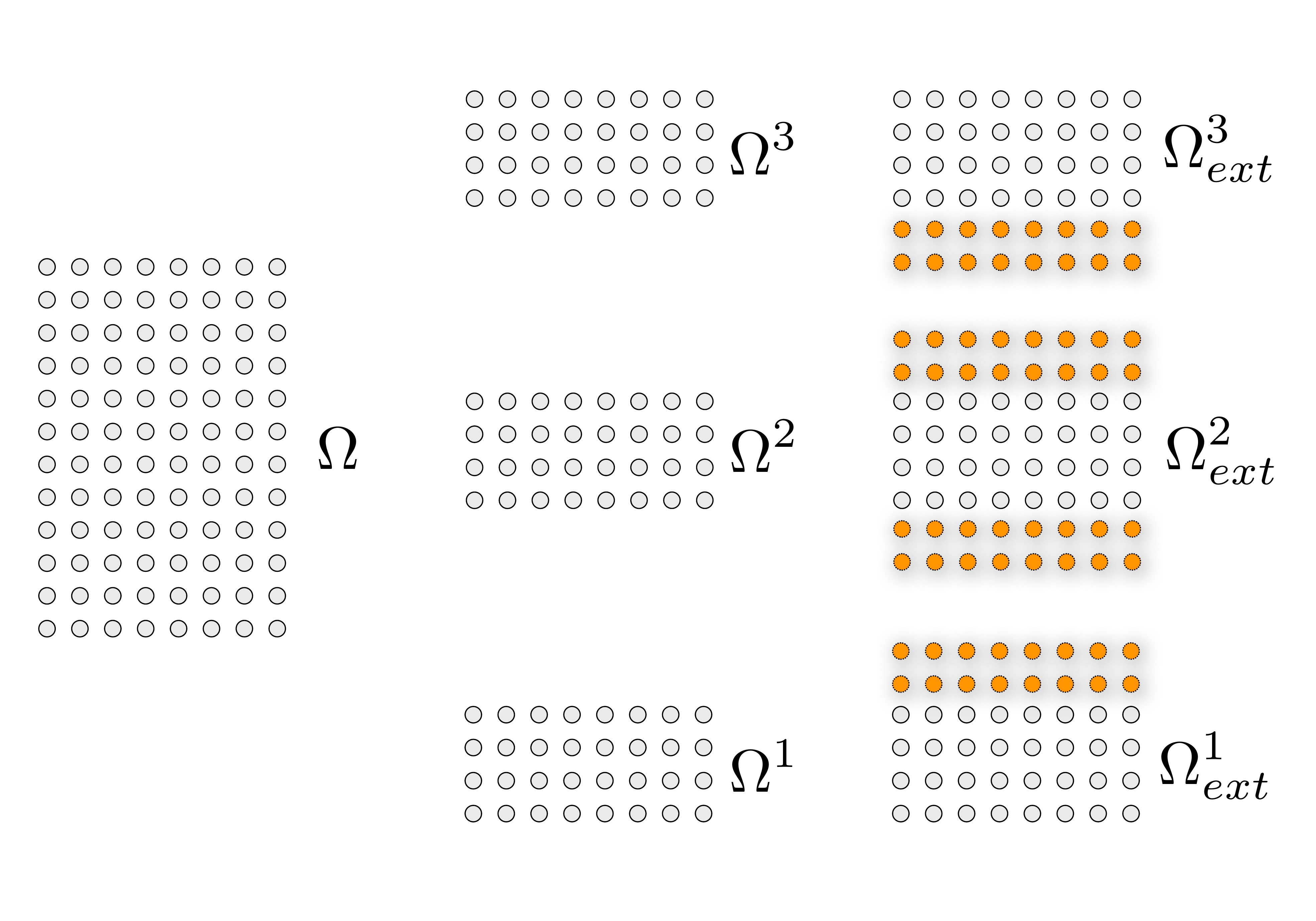}
\end{center}
\caption{Layered domain decomposition. The orange grid-points represent artificial grid points in the extended domain.}\label{fig:DDM_sketch}
\end{figure}

The domain $\Omega$ is partitioned in $L$ horizontal layers $\{ \Omega^{\ell} \}_{\ell =1}^L$ as depicted in Fig.~\ref{fig:DDM_sketch}. Following \cite{ZepedaDemanet:the_method_of_polarized_traces} we extend each subdomain, which we denote by $\{ \Omega^{\ell}_{ext} \}_{\ell =1}^L$.
Within each extended subdomain we assemble the discrete local problem
\begin{equation} \label{eq:local_discrete_sparse_pde}
	\C^{\ell} \w^{\ell} = \g^{\ell} = \chi_{\Omega^{\ell}} \g, \qquad \mbox{for } \ell = 1,..., L;
\end{equation}
where $ \chi_{\Omega^{\ell}}$ is the characteristic function of $\Omega^{\ell}$.

We impose $\C^{\ell}$ to satisfy the {\it compatibility condition}, i.e., $\C^{\ell}$ coincides exactly with $\C$ in the interior of $\Omega^{\ell}$. Within  $\Omega^{\ell}_{ext} - \Omega^{\ell}$, i.e., the extension of the subdomain, we do not require $\C^{\ell}$  to coincide with $\C$, instead we use the extra degrees of freedom to impose an analogue to an absorbing boundary condition. Given that $\C$ in the sparsified Lippmann-Schwinger system \eqref{eq:global_sparse_system} does not arise from a direct discretization of the Helmholtz equation, the boundary condition for subdomains needs a different treatment, which is presented in Section \ref{section:ABC}.

Fig.~\ref{fig:DDM_sketch} shows an example of such decomposition. At the left we have all the degrees of freedom, in the middle we have the domain decomposed in several layers, and at the right we have the extended domains. We represent in gray the grid points whose associated entries of each local matrix remain the same, and in orange the grid points that correspond to the extended degrees of freedom. Note that the entries of $\C^{\ell}$ associated to these extended grid points will be modified in order to impose an ABC.

We point out that the bidirectional preconditioner uses two of such decompositions: one in which the domain is decomposed in horizontal layers, such as in Fig.~\ref{fig:DDM_sketch}, and one in which the domain is decomposed in vertical layers. Alternating sweeping will be performed on these two settings.

\subsection{Algebraic decomposition}

As for any domain decomposition method, the objective is to solve the global problem \eqref{eq:global_sparse_system} by solving the local problems \eqref{eq:local_discrete_sparse_pde} iteratively, which requires transmission conditions at the interfaces between subdomains. In the case of high-frequency problems, due to possible pollution errors \cite{Sauter_Babuska:Is_the_Pollution_Effect_of_the_FEM_Avoidable_for_the_Helmholtz_Equation_Considering_High_Wave_Numbers}, such transmissions conditions need to be of higher order or algebraically exact with respect to the discrete global problem, in order to ensure a fast convergence. One advantage of the method of polarized traces is that such conditions can be seamlessly imposed at an algebraic level, bypassing the need to define optimized transmission conditions \cite{Gander:Optimized_Schwarz_Methods,Gander_Nataf:Optimized_Schwarz_Methods_without_Overlap_for_the_Helmholtz_Equation}.

In order to impose algebraically exact transmission conditions between layers, we use the blocks of $\C$ to define a discrete Green's representation formula. Such a formula can be deduced by imposing a discontinuous local solution, as it was performed in \cite{Stolk:An_improved_sweeping_domain_decomposition_preconditioner_for_the_Helmholtz_equation,Zepeda_Demanet:Nested_domain_decomposition_with_polarized_traces_for_the_2D_Helmholtz_equation}.

By construction we have that
\begin{equation} \scriptsize \label{eq:Helmholtz_matrix}
    \C =    \left [ \begin{array}{cccccc}
                    \C_{1,1}    &  \C_{1,2}     &           &                   &       \\
                    \C_{2,1}    &  \C_{2,2}     &  \C_{2,3} &                   &       \\
                                &     \ddots    &   \ddots  & \ddots            &    \\
                                &               &    \ddots & \ddots            &   \C_{n_z-1, n_z}  \\
                                &               &           &   \C_{n_z, n_z-1} &   \C_{n_z,n_z}
                    \end{array}
            \right ],
\end{equation}

 in which all the sub-matrices are sparse. Given that the mesh is uniform we use the following ordering of the unknowns
\begin{equation} \label{eq:global_ordering}
\v = (\v_1, \v_2, ..., \v_{n_z}),
\end{equation}
where $\v_j \in \CC^{n_x}$ is given by
\begin{equation} \label{eq:trace_ordering}
\v_j = (v_{1,j}, v_{2,j},..., v_{n_x,j}),
\end{equation}
which corresponds to $\v$ sampled at a constant $z$.

We write $\v^{\ell}$ for the wavefield defined locally at the $\ell$-th layer, i.e., $\v^{\ell}  = \chi_{\Omega^{\ell}}\v$, and $\v^{\ell}_k$ for the values at the local depth $z_k^{\ell}$ of $\v^{\ell}$.

In particular, $\v^{\ell}_1$ represents the wavefield sampled at the bottom of the layer, and  $\v^{\ell}_{n^{\ell}}$ is the wavefield at the top of the layer. Moreover,  $\v^{\ell}_0$ is the local wavefield sampled immediately below the bottom of the layer, in the extended region within the PML; $\v^{\ell}_{n^{\ell}+1}$  is the local wavefield sampled immediately above the top of the layer.

Finally, we use the notation $\delta_{k}$ to denote a discrete measure. For a generic vector $ \mathbf{p} \in \CC^{n_x}$,   $\delta_{k} \mathbf{p} \in \CC^{n_x\times n_z}$ is defined as follows:
\begin{equation}
 	(\delta_{k} \mathbf{p})_{i,j} = \left \{  \begin{array}{cl}
 									\mathbf{p}_i & \mbox{if } j = k,  \\
 									0        & \mbox{if } j\neq k.
 							\end{array} \right .
 \end{equation}
This notation will be used extensively to characterize the forcing terms necessary to transfer the information from one subdomain to its neighbors.

The correct Green's representation formula is given by the equivalent and more manageable expression\footnote{It is possible to recover the Green's representation formula by multiplying \eqref{eq:pde_GRF_real} by the inverse of $ \C^{\ell} $}
 \begin{align} \label{eq:pde_GRF_real}
    \C^{\ell} \w^{\ell} =  & - \delta_{n^{\ell}} \C^{\ell}_{n^{\ell},n^{\ell}+1} \v^{\ell}_{n^{\ell}+1}  +  \delta_{n^{\ell}+1}  \C^{\ell}_{n^{\ell}+1,n^{\ell}} \v^{\ell}_{n^{\ell}} \\
    & -  \delta_{1}  \C^{\ell}_{1,0} \v^{\ell}_{0}  + \delta_{0}  \C^{\ell}_{0,1} \v^{\ell}_{1} + \g^{\ell} . \nonumber
\end{align}

It can be shown that if the traces of the global $\v$, i.e. $\C \v = \f$, are known and used in \eqref{eq:pde_GRF_real}, then $ \w^{\ell}$ coincides exactly with $\v$ inside the layer, i.e., $ \w = \chi_{\Omega^{\ell}} \v$ (see Appendix C of \cite{ZepedaDemanet:the_method_of_polarized_traces} for a rigorous proof for matrices arising from finite differences). In  \cite{Zepeda_Demanet:Nested_domain_decomposition_with_polarized_traces_for_the_2D_Helmholtz_equation} the above formula is used to build an equivalent surface integral formulation in a matrix-free fashion, which is then solved iteratively. The convergence is fast provided that high-quality ABC are implemented, which is the topic of the next subsection.

\subsection{Absorbing boundary conditions} \label{section:ABC}
One of the most important elements of the method of polarized traces is the implementation of ABC at the interfaces between subdomains. ABC have successfully been used in domain decomposition for elliptic equations \cite{Engquist_Zhao:Absorbing_boundary_conditions_for_domain_decomposition}, and a simpler version of ABC were extensively used in the first domain decomposition method for the Helmholtz equation to ensure convergence \cite{Despres:domain_decomp}. High-quality absorbing conditions between subdomains are key to the fast convergence of iterative solvers for the high-frequency Helmholtz equation \cite{Chen_Xiang:a_source_transfer_ddm_for_helmholtz_equations_in_unbounded_domain,Cheng_Xiang:A_Source_Transfer_Domain_Decomposition_Method_For_Helmholtz_Equations_in_Unbounded_Domain_Part_II_Extensions,EngquistYing:Sweeping_PML,Leng:A_Fast_Propagation_Method_for_the_Helmholtz_equation,Leng_Ju:An_Overlapping_Domain_Decomposition_Preconditioner_for_the_Helmholtz_equation,Liu_Ying:Additive_Sweeping_Preconditioner_for_the_Helmholtz_Equation,Liu_Ying:Recursive_sweeping_preconditioner_for_the_3d_helmholtz_equation,CStolk_rapidily_converging_domain_decomposition,Stolk:An_improved_sweeping_domain_decomposition_preconditioner_for_the_Helmholtz_equation,GeuzaineVion:double_sweep}.

The high-quality ABC are needed in order to minimize the artificial reverberations at the interfaces between subdomains that arise from the truncation introduced by the domain decomposition. One of the objectives of the ABC, is to make the local Green's functions, which are defined by the local problem \eqref{eq:local_discrete_sparse_pde}, as close as possible to the global Green's function restricted to the  subdomain without scattering outside the subdomain. In other words, the local Green's functions should capture the local interactions within the subdomain, but should not contain pollution effects due to artificial reflections at the truncated boundary. Given that the medium outside the subdomain is not involved in the local system, the local Green's functions do not capture those long range interactions among subdomains; these interactions will be captured in an iterative manner during the sweeps.

In practice, without careful treatment at the truncated boundaries, the artificial reverberations can heavily pollute the local Green functions and hinder the convergence of the iterative method. For example, the adverse effect of the artificial reverberations can be overwhelming for the cases when Dirichlet or Neumann boundary conditions are imposed at the boundaries between subdomains \cite{GeuzaineVion:double_sweep,Gander:why_is_difficult_to_solve_helmholtz_problems_with_classical_iterative_methods}, causing the method to converge in $\cO(\omega)$ iterations, thus resulting in a super-linear complexity.

Given that the matrix $\C$ in our problem does not arise from a direct discretization of the Helmholtz equation, absorbing boundary conditions are not readily available. A natural way is to localize the perturbation, $m(\x)$, to the subdomain via a smooth cut-off. The localization produces a local discrete system in the form of \eqref{eq:local_discrete_sparse_pde} that satisfies the compatibility conditions.

Now, we proceed to give an explanation of how to define the local problems \eqref{eq:local_discrete_sparse_pde} using the localization to obtain an approximated ABC. Intuitively, the solution to \eqref{eq:local_discrete_sparse_pde} should provide an approximation of the wave scattered by the perturbation supported in $\Omega^{\ell}$. Thus, one naive approach would be to define the local problem as
\begin{equation}
	\cH^{\ell}u^{\ell} := u^{\ell} + \omega^2 G*(\chi_{\Omega^{\ell}} m u^{\ell}),
\end{equation}
and then use the sparsifying preconditioner to obtain a sparse local problem in the form of  \eqref{eq:local_discrete_sparse_pde}.

However, the local perturbation $\chi_{\Omega^{\ell}} m$ is discontinuous at the interfaces between subdomains, which will produce strong artificial reflections for waves impinging on the artificial boundary of $\Omega^{\ell}$. These strong reflections will pollute the local Green's function and result in a slow convergence. Fortunately, this adverse effect can be easily rectified by using a window, or cut-off, function that transitions the perturbation smoothly to zero.

Define a set of smooth cut-off functions $\{ \xi^{\ell}(\x) \}_{\ell = 1}^{L}$, such that
\begin{equation}
	\xi^{\ell}(\x) = \left \{  \begin{array}{cl}
							   1 &, \mbox{if } \x \in \Omega^{\ell};\\
							   0 &, \mbox{if } \x \notin \Omega^{\ell}_{ext}.
							   \end{array} \right .
\end{equation}
We can define the continuous counterpart of \eqref{eq:local_discrete_sparse_pde} as the problem
\begin{equation} \label{eq:local_Lippmann-Schwinger}
	\cH^{\ell}u^{\ell} := u^{\ell} + \omega^2 G*(\xi^{\ell} m u^{\ell}).
\end{equation}
which is posed for $\x \in \Omega$.
We recall that the matrix $\mathbf{D}$, is a diagonal matrix that encodes the multiplication by $m(\x_{i,j})$; we define $\mathbf{D}^{\ell}$ as the diagonal matrix that represents the multiplication by $\xi^{\ell}(\x)m(\x)$. We can define the discrete local Lippmann-Schwinger equation as,
\begin{equation} \label{eq:whole_discrete_local_Lippmann-Schwinger}
\H^{\ell}_{global} = \I + \omega^2 \mathbf{G} \cdot \mathbf{D}^{\ell}.
\end{equation}
The global subscript indicates that $\H^{\ell}_{global}$ is a $N \times N$ matrix. In this case we can multiply \eqref{eq:whole_discrete_local_Lippmann-Schwinger} by $\mathbf{A}$, and we would obtain a local system that satisfies the compatibility condition of \eqref{eq:local_discrete_sparse_pde}. However, this would be prohibitively expensive given that would need to solve a $N \times N$ system for each subdomain.

The key observation is that only a small portion of  $\H^{\ell}_{global}$  in \eqref{eq:whole_discrete_local_Lippmann-Schwinger}  deviates from an identity matrix. This is caused by the small support of $\xi^{\ell} m$, in which $ \mbox{supp}(\xi^{\ell} m )\subset \Omega^{\ell}_{ext}$. Thus, we can re-write \eqref{eq:local_Lippmann-Schwinger} in $\Omega^{\ell}_{ext}$ instead of $\Omega$, and discretize it resulting in
\begin{equation} \label{eq:discrete_local_Lippmann-Schwinger}
\H^{\ell} = \I + \omega^2 \mathbf{G} \cdot \mathbf{D}^{\ell},
\end{equation}
which is a quasi-1D (see \cite{EngquistYing:Sweeping_PML}) $q n \times q n $ system, where $q$ is the thickness, in grid points, of each extended subdomain. In the implementation of the method, the cut-off is constant in the direction tangential to the boundary, $x$ direction in our layered setting, and it is a third order spline in the direction normal to the boundary, $z$ direction in our setting.

In theory, we could apply the sparsifying preconditioner to \eqref{eq:discrete_local_Lippmann-Schwinger}; however, given that the mesh is also truncated, the optimization via a SVD can provide different results, hindering the compatibility condition. Instead, we multiply \eqref{eq:discrete_local_Lippmann-Schwinger} by the corresponding block of $\mathbf{A}$, but with a special care at the boundaries of $\Omega^{\ell}_{ext}$, where we use the stencil generated at the boundary of $\Omega$, which provides a satisfying ABC \cite{Ying:Sparsifying_Preconditioner_for_the_Lippmann--Schwinger_Equation} for the background wave-speed at the outer boundary of the extended subdomain.

\begin{figure}[h]
\begin{center}
\includegraphics[trim= 172mm 55mm 55mm 50mm, angle =0,clip, width=13cm]{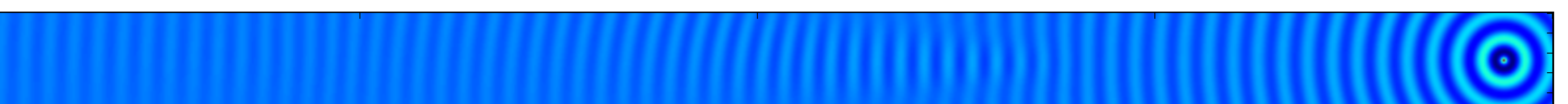} \\
\vspace{0.2cm}
\includegraphics[trim= 172mm 55mm 55mm 50mm, angle =0,clip, width=13cm]{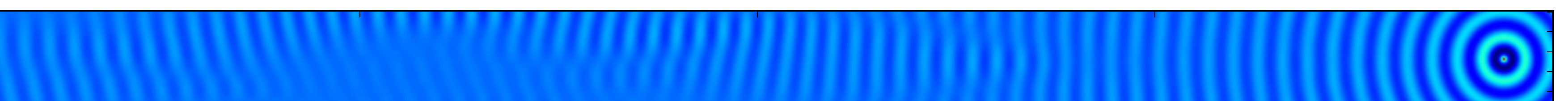}
\end{center}
\caption{Example of a global (top) and local (bottom) Green's function for a positive perturbation, we can observe that there is some interference in the local Green's function given by the spurious grazing wave that is concentrated at the top of the slab.}\label{fig:Green_functions+}
\end{figure}

\begin{figure}[h]
\begin{center}
\includegraphics[trim= 172mm 55mm 55mm 50mm, angle =0,clip, width=13cm]{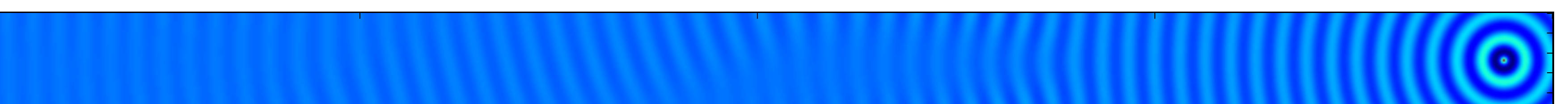} \\
\vspace{0.2cm}
\includegraphics[trim= 172mm 55mm 55mm 50mm, angle =0,clip, width=13cm]{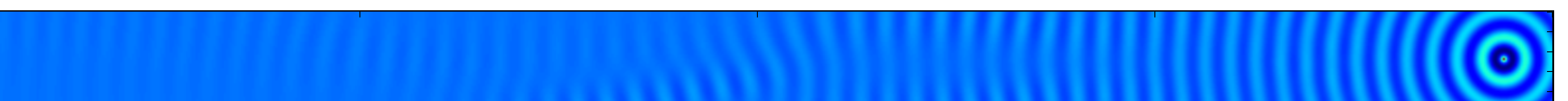}
\end{center}
\caption{Example of a global (top) and local (bottom) Green's function for a negative perturbation, we can observe some differences between the Green's functions, but we do not observe the kind of interference, nor the concentration at the boundaries, as observed in Fig.~\ref{fig:Green_functions+}.}
\label{fig:Green_functions-}
\end{figure}

Unfortunately, the smooth cut-off used to localize the perturbation is not equivalent to an ABC, which would approximate the Dirichelet to Neumann (DtN) map at the boundary; such that it can absorb and dampen waves propagating out of the boundary in all directions with little reflections. The smooth cut-off in the normal direction to the boundary minimally affects waves propagating out of the domain close to the normal direction. However, for waves propagating  near and almost tangential to the boundary, a smooth cut-off may not be good enough depending on the sign of the perturbation as illustrated in Figures \ref{fig:Green_functions+} and \ref{fig:Green_functions-}. If the perturbation is negative, the waves propagate faster inside the subdomain than in the cut-off region.
If a wave propagates almost tangentially to the boundary to the cut-off region it will bend outwards, and it will be ultimately absorbed by the numerical ABC situated at the artificial boundary of the extended subdomain. On the other hand, if the perturbation is positive, then the waves propagate slower inside the subdomain than in the background media. Under these circumstances, the subdomain with a smooth cut-off to the background media creates an artificial wave-guide effect at the subdomain boundary. Waves propagating almost normal to the boundaries  will still cross the smooth cut-off and get absorbed by the numerical ABC at the artificial boundary of the extended subdomain. However, waves propagating almost transversally along the subdomain boundary will bend inwards due to the artificial cut-off, producing spurious grazing waves that pollute the local solution within the subdomain. This last problem can be attenuated to some extent by introducing a $\cO(\omega)$ complex shift within the cut-off region, thus dampening the grazing waves. We point out that we are only introducing a complex shift in the extended domain; it would be possible to introduce a complex shift in $\Omega$ as in \cite{Graham:Domain_Decomposition_preconditioning_for_high-frequency_Helmholtz_problems_using_absorption}, however, that would violate the compatibility condition.

% We may want to put some examples and some pictures in here.

Although the complex shift helps to reduce the number of iterations for the original method of polarized traces from $\cO(\omega)$ to $\cO(\sqrt{\omega})$, which is better than standard iterative methods, it is not optimal. In order to obtain the lower rate of $\cO(\log{\omega})$ iterations to converge, we introduce in the next section an effective bi-directional alternating sweeping preconditioner to handle those grazing waves due to either the incoming waves or the smooth cut-off at the artificial boundaries between two subdomains in one directional sweeping.

\subsection{One-directional algorithm} \label{section:onedirectional_algorithm}

In this section we discuss the preconditioner for the sparse problem based on the method of polarized traces for one-directional sweep and we provide the algorithm in pseudo-code with its physical interpretation.

The main intuition behind the method of polarized traces is the directionality of the waves, which are imposed by the polarizing conditions (see Section 3 of \cite{ZepedaDemanet:the_method_of_polarized_traces}) and the absorbing boundary conditions. For example, the wavefield generated by a local source, $\g^{\ell}$
\begin{equation}
\C^{\ell} \v^{\ell} = \g^{\ell},
\end{equation}
irradiates from the interior of $\Omega^{\ell}$ to the exterior of it. If we sample $\v^{\ell}$ at the bottom of $\Omega^{\ell}$, we can physically interpret it as the trace of a wavefield propagating downwards.

Following this physical interpretation, the wavefield propagating downwards from $\Omega^{\ell}$ propagates inwards $\Omega^{\ell-1}$. By using the Green's representation formula we can propagate the wavefield from $\Omega^{\ell}$ to $\Omega^{\ell-1}$ by artificial sources located at the top of $\Omega^{\ell-1}$ following
\begin{align} \label{eq:pde_GRF_down}
    \C^{\ell-1} \v^{\ell-1} =  & - \delta_{n^{\ell-1}} \C^{\ell-1}_{n^{\ell-1},n^{\ell-1}+1} \v^{\ell}_{1}  +  \delta_{n^{\ell-1}+1}  \C^{\ell-1}_{n^{\ell-1}+1,n^{\ell-1}}  \v^{\ell}_{0}  + \g^{\ell-1},
\end{align}
which is an alternate expression for the incomplete Green's representation formula introduced in \cite{ZepedaDemanet:the_method_of_polarized_traces}. %For a more complete exposition please see Section 1.1 in \cite{ZepedaDemanet:the_method_of_polarized_traces}.

We can observe that $\v^{\ell-1} $ contains the energy radiating from $\g^{\ell-1}$ and $\g^{\ell}$. If we sample $\v^{\ell-1}$ at the bottom of $\Omega^{\ell-1}$, it can be interpreted as the downwards wavefield consisting of the downwards wavefield generated by $\g^{\ell-1}$ and the wavefield generated by $\g^{\ell}$ that propagated through $\Omega^{\ell-1}$.

We can continue the propagation by adding the local wavefields in a manner that is usually called a multiplicative Schwarz iteration (see \cite{Gander_Nataf:Optimized_Schwarz_Methods_without_Overlap_for_the_Helmholtz_Equation}). Analogously, the same procedure can be performed for the wavefields propagating upwards. The application of both sweeps, downwards and upwards can be easily understood as a block Jacobi iteration as presented in \cite{ZepedaDemanet:the_method_of_polarized_traces}. Furthermore, the polarizing conditions can be further exploited to reduce the number of iterations by using a Gauss-Seidel approach described by Alg.~\ref{alg:GS_preconditioner}, which in practical terms reduces the number of iterations by a factor of two (see \cite{Zepeda_Demanet:Nested_domain_decomposition_with_polarized_traces_for_the_2D_Helmholtz_equation}).

All the intuitive arguments provided in the current exposition can be made rigorous from an algebraic point of view, we skip such explanation and redirect the interested readers to \cite{ZepedaDemanet:the_method_of_polarized_traces}.

\begin{algorithm} \small Gauss-Seidel Preconditioner \label{alg:GS_preconditioner}
    \begin{algorithmic}[1]
        \Function{ $\v$ = Preconditioner}{ $\g$ }
            \For{  $\ell = L:-1:1$ } \Comment{Downwards Sweep}
                \State $ \g^{\ell} = \g \chi_{\Omega^{\ell}} $ \Comment{Partition of the Source}
                \If{ $\ell \neq 1$ }
                %\State $ \g^{\ell} = \g^{\ell} - \delta_{1}  \C^{\ell}_{1,0} \u_{0}  + \delta_{0}  \C^{\ell}_{0,1} \u_{1} $
                \State $ \g^{\ell} +=  - \delta_{n^{\ell}} \C^{\ell}_{n^{\ell},n^{\ell}+1} \mathbf{w}^{\ell-1}_{1}  +  \delta_{n^{\ell}+1}  \C^{\ell}_{n^{\ell}+1,n^{\ell}}  \mathbf{w}^{\ell-1}_{0}  $
                \EndIf
                \State $\mathbf{w}^{\ell}  = (\mathbf{C}^{\ell})^{-1}\g^{\ell}$ \Comment{Local Solve}
            \EndFor

            \For{  $\ell = 1: L$ } \Comment{Upwards Sweep}
                %\State $ \g^{\ell} = \f \chi_{\Omega^{\ell}} $ \Comment{Partition of the Source}
                \If{ $\ell \neq L$ }
                %\State $ \g^{\ell} = \g^{\ell} - \delta_{1}  \C^{\ell}_{1,0} \u_{0}  + \delta_{0}  \C^{\ell}_{0,1} \u_{1} $
                \State $ \g^{\ell} += - \delta_{1} \C^{\ell}_{1,0} \v^{\ell+1}_{n^{\ell}}  + \delta_{0}  \C^{\ell}_{0,1} \v^{\ell+1}_{n^{\ell}+1} $
                \EndIf

                %\If{ $\ell \neq 1$ }
                %\State $ \g^{\ell} = \g^{\ell} - \delta_{1}  \C^{\ell}_{1,0} \v_{0}  + \delta_{0}  \C^{\ell}_{0,1} \v_{1} $
                %\State $ \g^{\ell} +=  - \delta_{n^{\ell}} \C^{\ell}_{n^{\ell},n^{\ell}+1} \mathbf{w}^{\ell-1}_{1}  +  \delta_{n^{\ell}+1}  \C^{\ell}_{n^{\ell}+1,n^{\ell}}  \mathbf{w}^{\ell-1}_{0}  $
                %\EndIf
                \State $\v^{\ell}  = (\mathbf{C}^{\ell})^{-1}\g^{\ell}$  \Comment{Local Solve}
            \EndFor

            \State $\v =  \left (\v^{1}, \v^{2}, \v^{3}, ..., \v^{L-1}, \v^{L}  \right)$ \Comment{Concatenation}
        \EndFunction
    \end{algorithmic}
\end{algorithm}

\subsection{Bi-directional algorithm} \label{section:bidirectional_algorithm}

In this Section we first discuss briefly the shortcomings of Alg.~\ref{alg:GS_preconditioner} when applied to the sparsified Lippmann-Schwinger equation and then introduce the bi-directional preconditioner in order to achieve the optimal complexity.

The algorithm as presented in Alg.~\ref{alg:GS_preconditioner}, may not work optimally in our setting for the Lippmann-Schwinger equation, in which the perturbation is localized using a smooth cut-off. The main issue is that waves propagating nearly tangential to the artificial boundaries between subdomains are perturbed due to the cut-off and depending on the perturbation they can be bended inwards and thus generate spurious grazing waves. Even though the spurious grazing waves are dampened by the introduction of a complex shift in the cut-off region to some extent, they can still cause significant deviations of the local Green's functions with respect to the global ones. Thus, they hinder the fast convergence of the iterations, especially when the perturbation is positive, i.e., the wave-guide situation, as explained in Section \ref{section:ABC} and shown in Figures \ref{fig:Green_functions+} and \ref{fig:Green_functions-}.

Moreover, a closer inspection to Alg.~\ref{alg:GS_preconditioner} reveals that within the algorithmic pipeline, the residual given by $\mathbf{r}  = \mathbf{C} \v_{precond} - \g $, where $ \v_{precond} = Preconditioner(\g)  $ are concentrated at the boundaries between subdomains. In particular, $\mathbf{r}$ can be interpreted in the continuous case as a measure supported at the boundaries between subdomains with most singular term of the form $\delta'(\x)$. Then, as explained in Page 243 of \cite{CStolk_rapidily_converging_domain_decomposition}, the restriction of such measure to a subdomain (line 3 in Alg.~\ref{alg:GS_preconditioner}) is not well defined; such ambiguity ultimately results in a slow convergence rate. One way to solve the residual problem is presented in \cite{CStolk_rapidily_converging_domain_decomposition,Stolk:An_improved_sweeping_domain_decomposition_preconditioner_for_the_Helmholtz_equation,Zepeda_Hewett_Demanet:Preconditioning_the_2D_Helmholtz_equation_with_polarized_traces}, where an overlapping domain decomposition is used to ensure that at each iteration the residual lies at the interior of a subdomain.

% Moreover, within the algorithm pipeline, the residuals are concentrated at the boundaries after the application of the one directional sweeping preconditioner. In the next iteration, during the decomposition of the source, the residuals should be split between two subdomains. However, there is not a simple way to decompose the grazing waves near the boundaries correctly like the polarized traces for impinging waves, thus resulting in an ambiguity, which in return results in a slow convergence rate. One way to solve the residual problem is presented in \cite{CStolk_rapidily_converging_domain_decomposition,Stolk:An_improved_sweeping_domain_decomposition_preconditioner_for_the_Helmholtz_equation,Zepeda_Hewett_Demanet:Preconditioning_the_2D_Helmholtz_equation_with_polarized_traces}, where a overlapping domain decomposition is used to ensure that at each iteration the residual lies at the interior of a subdomain.

In this work, we solve both problems simultaneously by introducing a simple alternating bi-directional preconditioner, that consists on applying an up-down sweep performed by Alg.~\ref{alg:GS_preconditioner} followed by an left-right pass sweep, performed by Alg.~\ref{alg:GS_preconditioner} but applied to the transpose of $\C$, which is denoted $\mbox{Preconditioner}^T$ in Alg.~\ref{alg:birirectional_preconditioner}. The introduction of the bi-directional preconditioner in Alg.~\ref{alg:birirectional_preconditioner} handles effectively the spurious grazing waves in one sweep by performing another sweep in the orthogonal direction. In addition, the application of the preconditioner in the orthogonal direction, after a step of iterative refinement confines the residual to the interior of the subdomains, thus eliminating the ambiguity and resulting in a fast convergence.

In particular, as it will be shown in the numerical tests, the bi-directional preconditioner is more efficient and less sensitive to the geometry, direction of incoming waves and sign of the perturbation, albeit with a higher memory footprint.
The sweep in the orthogonal direction can easily be implemented as preconditioning the transpose of the matrix $\mathbf{C}$ using Alg.~\ref{alg:GS_preconditioner}.
\begin{algorithm} \small Bi-directional Preconditioner \label{alg:birirectional_preconditioner}
    \begin{algorithmic}[1]
        \Function{ $\v$ = BiDirectionalPreconditioner}{ $\g$ }
      	\State $\u_1$ = Preconditioner $(\g)$   \Comment{Up-Down pass}
        \State $\mathbf{e} = \g -  \mathbf{C}\u_1$  \Comment{Iterative refinement}
        \State $\u_2$ = Preconditioner$^{T}$ $(\mathbf{e})$  \Comment{Left-Right pass}
        \State $\v = \u_1 + \u_2$
        \EndFunction
    \end{algorithmic}
\end{algorithm}

\section{Two-level Preconditioner} \label{section:preconditioner}
We briefly summarize the operations within the two levels of the preconditioner presented in this paper as follows:

\begin{itemize}
\item the inner level, in which $u_{approx}$ in \eqref{eq:preconditioner} is computed approximately using the bi-directional preconditioner Alg.~\ref{alg:birirectional_preconditioner};

\item and, the outer level, in which \eqref{eq:discrete_Lipmann-Schwinger} is solved iteratively using GMRES preconditioned with Alg.~\ref{alg:Lipmann-Schwinger_preconditioner}.
\end{itemize}

\begin{algorithm} \small Lippmann-Schwinger Preconditioner \label{alg:Lipmann-Schwinger_preconditioner}
    \begin{algorithmic}[1]
        \Function{ $\u$ = Lippmann-Schwinger Preconditioner}{ $\f$ }
          \State $ u = \mbox{GMRES}( \C ,\mathbf{A} \f, BiDirectionalPreconditioner, tol = 10^{-3})   $
        \EndFunction
    \end{algorithmic}
\end{algorithm}

There is a balance to the choice of the tolerance for the inner GMRES iteration (in Alg.~\ref{alg:Lipmann-Schwinger_preconditioner}): if the tolerance is chosen too small, then the constants are large; on the other hand, if the tolerance is too big, the constants will decrease but the number of outer iterations will grow with the frequency, thus increasing the asymptotic cost. The best choice in our case, is to set the inner loop tolerance to $10^{-3}$. This choice ensures that the amount of outer iterations will be the same as using a direct multi-frontal method and produces reasonable constants.

As mentioned in the introduction, the algorithm presented in this paper has two phases: an off-line phase performed only once, and an on-line phase performed for each right-hand-side, or incoming wave. The complexity of the solver is summarized in Table \ref{table:complexity}; in the sequel we provide further details for the complexity for each phase, in which we suppose the number of subdomains, L, is $\cO(n)$ and that $q = \cO(1)$.

\begin{table}  \small
\begin{center}
\begin{tabular}{|c|c|}
\hline
Stage    & Complexity  \\
\hline
Off-line & $\cO\left(  N\right)$   \\
\hline
On-line  &  $\cO \left(n_{\mbox{iter}}^{\mbox{sp}} \left( N\log{N} + N \cdot n_{\mbox{iter}} \right)  \right)$  \\
\hline
\end{tabular}
\end{center}
\caption{Complexity for each state of the solver, supposing the $q$ ,the thickness of each extended local subdomain, is $\cO(1)$.}
\label{table:complexity}
\end{table}

\subsection{Off-line complexity}

For the off-line or setup phase, we need to build $\C$ for which we need to compute $\mathbf{A}$. As explained in \cite{Ying:Sparsifying_Preconditioner_for_the_Lippmann--Schwinger_Equation} $\mathbf{A}$ can be constructed by solving only a hand-full of optimization problems. This property arises mainly by the translation invariance of the convolution kernel, meaning that we can reuse the same stencil for many grid-points. Solving the minimization problem needs a SVD of a $9 \times N$ matrix, which can be performed in $\cO(N)$ complexity. To create $\C$ we need to multiply $\mathbf{A}$ times $\H$; exploiting the sparsity of $\mathbf{A}$ and performing the multiplication in the correct order, the product can be carried in $\cO(N)$ operations.

We can build all the local matrices $\C^{\ell}$ in $\cO(N)$ complexity by reusing $\mathbf{A}$; we build $L= \cO(n)$ local matrices, and each $\C^{\ell}$ has $\cO(n)$ non-zero elements. Finally, the factorization of all the local problems can be performed in $\cO(N)$ complexity. Each local matrix is a circulant matrix and their factorization can be performed in $\cO(n)$ time, and we have $L =\cO(n)$ systems to factorize, which results in $\cO(N)$ complexity. We point out that assembly and factorization of the local matrices is an embarrassingly parallel operation which can be effortlessly parallelized.

\subsection{On-line complexity}

For the outer GMRES we need to apply $\H$, which consists of multiplications by sparse matrices and a convolution that can be performed fast via a FFT (see section 2.3 of \cite{Duan-Rohklin:High-order_quadratures_for_the_solution_of_scattering_problems_in_two_dimensions}). The cost of applying $\H$ is dominated by the application of the FFT, resulting on an overall $\cO(N \log{N})$ complexity. For Alg.~\ref{alg:GS_preconditioner} each local solve can be performed in $\cO(n)$ time resulting in a complexity $\cO(N)$. We have an overall complexity of $\cO(N)$ for each application of Alg.~\ref{alg:GS_preconditioner}.

Finally, the overall on-line complexity depends on the number of outer and inner iterations to converge. From \cite{Ying:Sparsifying_Preconditioner_for_the_Lippmann--Schwinger_Equation} the number of outer iterations is almost independent of the frequency, $n_{\mbox{iter}}^{\mbox{sp}}  = \cO(1)$, and we provide in the next section numerical evidence that under the assumptions mentioned in the introduction, the number of inner iterations grows logarithmically with the frequency, $ n_{\mbox{iter}} = \cO(\log{\omega})$, and given that we are in the high-frequency regime $ n_{\mbox{iter}} = \cO(\log{N})$, resulting in the overall complexity mentioned in the introduction.

% \subsection{Complexity with respect to $\omega$}(???)

% In this section we could briefly explain about the complexity of the algorithm with respect to the frequency. There is not pollution error and we are only exposed to the accuracy of the quadratures

\section{Numerical Experiments} \label{section:numerical_experiments}

The numerical experiments were implemented in Julia v0.4.3 linked against Intel MKL 11.2, and were performed in a server with dual socket Intel Xeon CPU E5-2670, and 384 Gb of RAM. We used the FFT library FFTW3 \cite{FFTW:Johnson_Frigo} for the application of $\H$. For the local solves in the inner iteration we used the PARDISO \cite{Schenk:Solving_unsymmetric_sparse_systems_of_linear_equations_with_PARDISO} and UMFPACK \cite{Davis:UMFPACK} solvers.

One of the main advantages of the method of polarized traces lies on its modularity. We used the state-of-the-art library FFTW3 for the application of $\H$ but it is possible to further accelerate the application of $\mathbf{H}$ by using GPU-accelerated libraries (such as cuFFT \cite{cuda}), or a distributed-memory FFT libraries. Within the main bottleneck of the algorithm, which is the application of the preconditioner, it is possible to effortlessly exchange the linear solver. In our experiments we used PARDISO \cite{Schenk:Solving_unsymmetric_sparse_systems_of_linear_equations_with_PARDISO} and UMFPACK \cite{Davis:UMFPACK}, but it is easy to exchange them by solvers using compressed linear algebra \cite{Rouet_Li_Ghysels:A_distributed-memory_package_for_dense_Hierarchically_Semi-Separable_matrix_computations_using_randomization,Sivaram_Darve:HODLR,Xia:multifrontal,Gillman_Barnett_Martinsson:A_spectrally_accurate_solution_technique_for_frequency_domain_scattering_problems_with_variable_media}, or by massively parallel solvers such as MUMPS \cite{Amestoy_Duff:MUMPS} (or its compressed linear algebra variant \cite{Amestoy_Weisbecker:compressed_MUMPS}) or SuperLU-DIST \cite{Demmel_Li:superlu}, among many others.

In this section we corroborate the complexity claims for our method; in particular, present 3 different kinds of numerical experiments:
\begin{enumerate}
\item we study the asymptotic cost of solving sparsified systems from discretizations of the Lippmann-Schwinger equation using the bi-directional preconditioner;
\item we solve a wave scattering example in three different medium profiles shown in Fig.~\ref{fig:scattering} using the full two level preconditioner,
\item  we solve an example arising from reflectometry of plasmas confined in a fusion reactor \cite{Weitzner:Lower_hybrid_waves_in_the_cold_plasma_model}.
% in which a wave propagates until it meets a zone in which the equation becomes strictly elliptic, thus reflecting back.
\end{enumerate}

\subsection{Sparsified system}
\begin{figure}[h]
\begin{center}
\includegraphics[trim= 19mm 10mm 23mm 15mm, angle =0,clip, width=6.3cm]{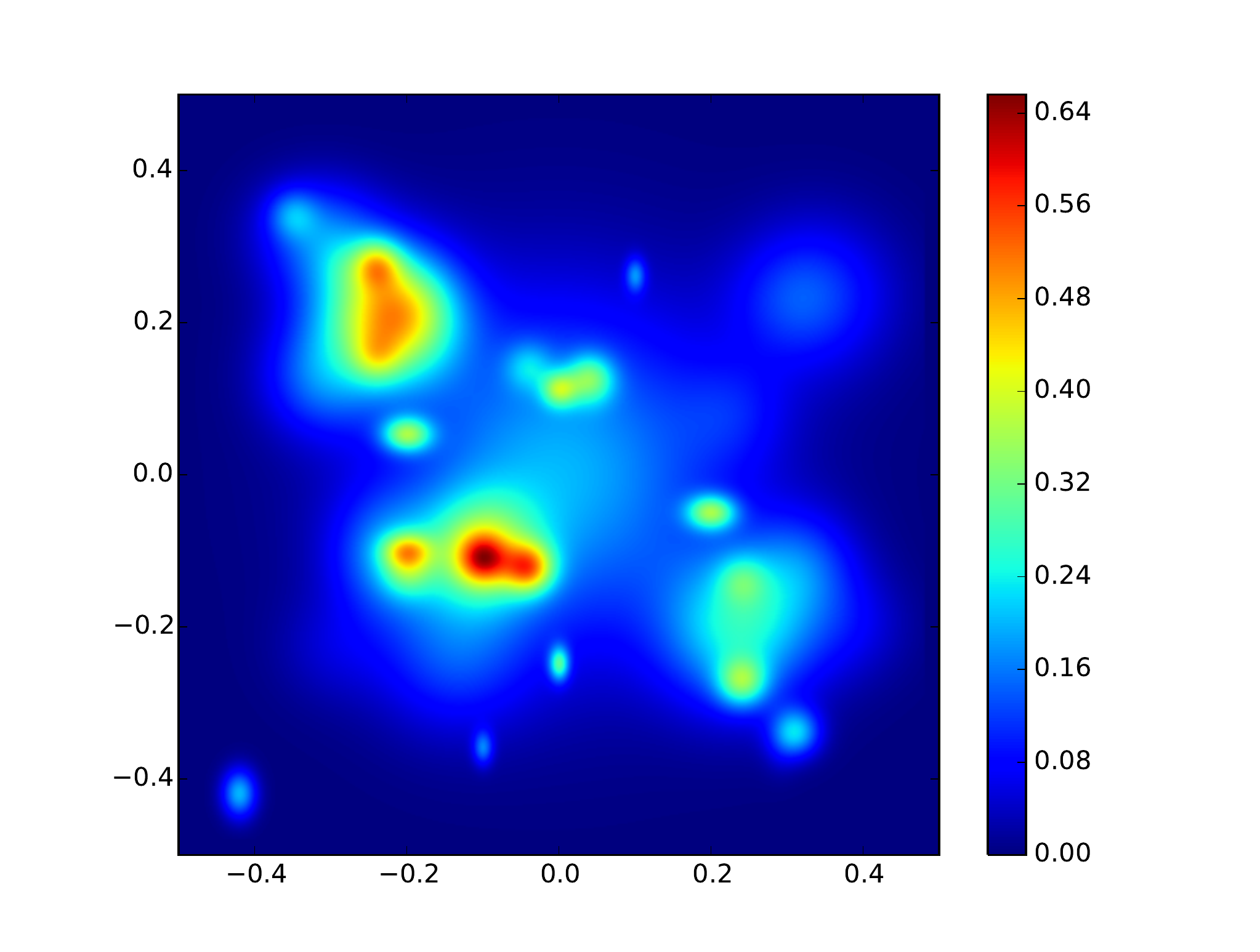}
\includegraphics[trim= 19mm 10mm 23mm 15mm, angle =0,clip, width=6.3cm]{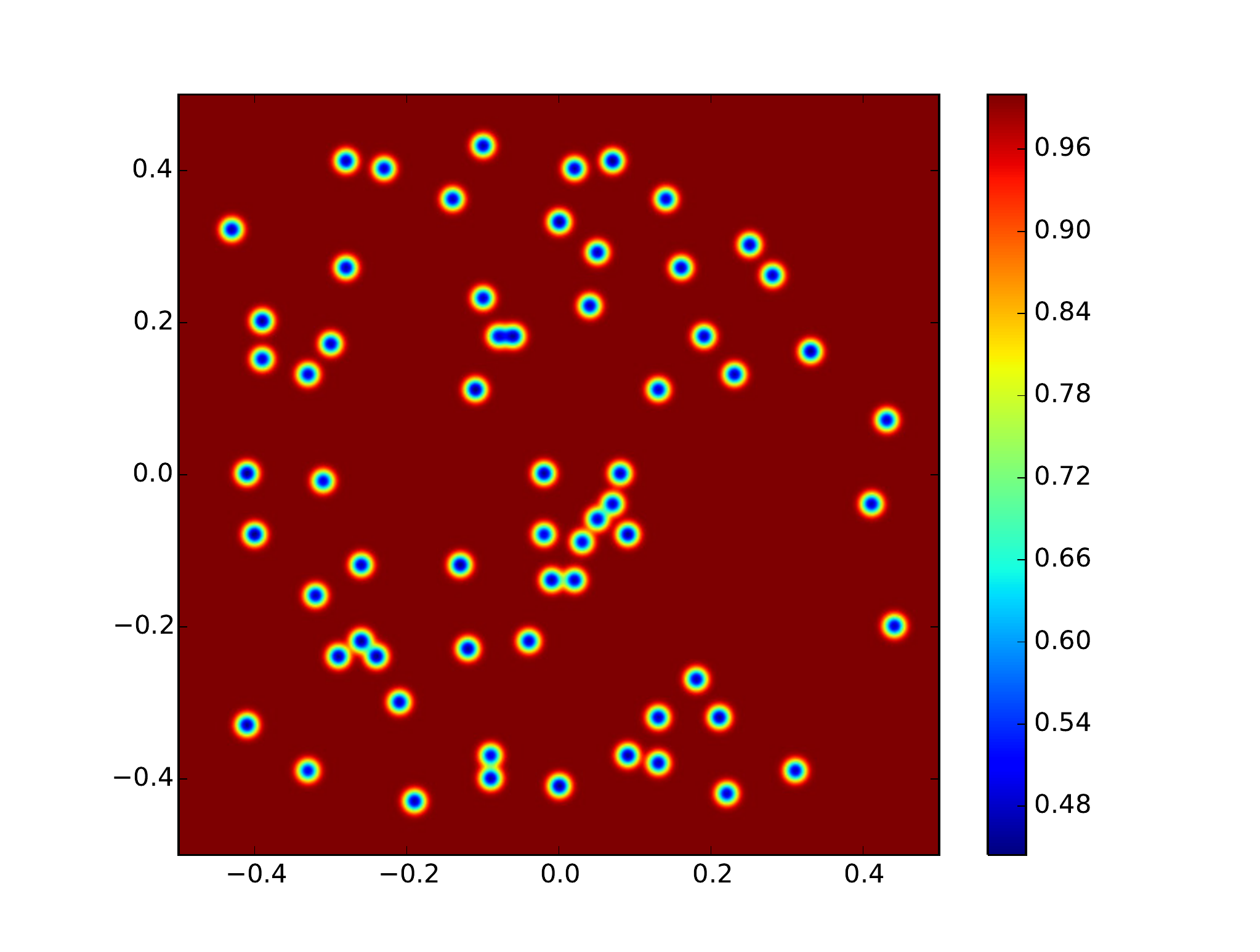}

\end{center}
\caption{(left) Perturbation used for the wave scattering. (right) Wave-speed for 64 randomly placed Gaussian bumps}\label{fig:scattering}
\end{figure}

First, we assemble two linear systems by sparsifying the Lippmann-Schwinger equation for the perturbation given in Fig.~\ref{fig:scattering} ({\it left}) and demonstrate the performance of the original method of polarized traces, and the new bi-directional preconditioner. The two sparsified Lippmann-Schwinger systems \eqref{eq:global_sparse_system} correspond to the same perturbation shown in Fig.~\ref{fig:scattering} ({\it left}) but with a different sign, to which we  refer to as the positive and negative perturbation. We solve the resulting sparsified systems with GMRES (tolerance of $10^{-6}$) preconditioned with the method of polarized traces \cite{Zepeda_Demanet:Nested_domain_decomposition_with_polarized_traces_for_the_2D_Helmholtz_equation} and the bi-directional preconditioner in Alg.~\ref{alg:birirectional_preconditioner}. We solve \eqref{eq:global_sparse_system} for $64$ different right-hand sides representing $64$ different incident waves, and we compute an average time. We repeat the process for a set of increasing frequencies, such that we have $10$ points per wavelength in the background medium. Moreover, the domain is decomposed in $L = n/50$ subdomains, such that each subdomain is roughly $5$ wavelengths thick, and each subdomain is extended by 10 grid points in both directions (see Fig.~\ref{fig:DDM_sketch}). The average execution times are presented in Table \ref{table:numerical_experiment_bi_directional}.

\begin{table} \small
\begin{center}
\begin{tabular}{|c|c|c|c|c||c|c| }
\hline
  N    & $\omega$ & $L$ & $T^{+}_{PT}$ [s] & $T^{+}_{BD}$ [s] & $T^{-}_{PT}$ [s] & $T^{-}_{BD}$ \\
\hline
$200  \times 200  $ & $200$  &  $ 4  $ &  $  0.5$ $\mathbf{(3.6)}$ & $ 0.2$ $\mathbf{(2.7)}$ & $ 0.5$ $\mathbf{(3.7)}$ & $ 0.2$  $\mathbf{(2.9)}$ \\
$400  \times 400  $ & $400$  &  $ 8  $ &  $  2.1$ $\mathbf{(4.2)}$ & $ 0.7$ $\mathbf{(3.1)}$ & $ 2.2$ $\mathbf{(4.3)}$ & $ 0.6$  $\mathbf{(3.0)}$ \\
$500  \times 500  $ & $500$  &  $ 10 $ &  $  3.5$ $\mathbf{(4.6)}$ & $ 1.1$ $\mathbf{(3.2)}$ & $ 3.3$ $\mathbf{(4.5)}$ & $ 1.0$  $\mathbf{(3.0)}$ \\
$800  \times 800  $ & $800$  &  $ 16 $ &  $ 11.5$ $\mathbf{(5.7)}$ & $ 3.2$ $\mathbf{(3.8)}$ & $8.62$ $\mathbf{(4.8)}$ & $ 2.7$  $\mathbf{(3.1)}$ \\
$1000 \times 1000 $ & $1000$ &  $ 20 $ &  $ 20.4$ $\mathbf{(6.6)}$ & $ 5.1$ $\mathbf{(4.1)}$ & $14.0$ $\mathbf{(5.2)}$ & $ 4.2$  $\mathbf{(3.1)}$ \\
$1250 \times 1250 $ & $1250$ &  $ 25 $ &  $ 30.5$ $\mathbf{(7.9)}$ & $ 8.7$ $\mathbf{(4.3)}$ & $22.5$ $\mathbf{(5.4)}$ & $ 7.3$  $\mathbf{(3.3)}$ \\
$1400 \times 1400 $ & $1400$ &  $ 28 $ &  $ 43.5$ $\mathbf{(8.9)}$ & $11.1$ $\mathbf{(4.7)}$ & $29.1$ $\mathbf{(5.6)}$ & $ 8.7$  $\mathbf{(3.4)}$ \\
$1600 \times 1600 $ & $1600$ &  $ 32 $ &  $ 62.8$ $\mathbf{(10.1)}$ & $15.8$ $\mathbf{(4.7)}$ & $41.8$ $\mathbf{(5.8)}$ & $12.6$ $\mathbf{(3.4)}$ \\
$2000 \times 2000 $ & $2000$ &  $ 40 $ &  $120.5$ $\mathbf{(12.6)}$ & $25.0$ $\mathbf{(5.3)}$ & $63.8$ $\mathbf{(6.1)}$ & $18.1$ $\mathbf{(3.5)}$ \\
\hline
\end{tabular}
\end{center}
\caption{Average number of iterations (in parenthesis) and average execution time to solve the sparsified system arising from the Lipmann-Schwinger equation using the positive (+) and negative (-) perturbation, using either the method of polarized traces (PT), or the new bi-directional preconditioner(BD). We use GMRES with a tolerance of $10^{-6}$, in this case $N = n^2$ is the number of degrees of freedom, $\omega$ is the frequency and $L$ is the number of subdomains.} \label{table:numerical_experiment_bi_directional}

\end{table}

Table \ref{table:numerical_experiment_bi_directional} clearly depicts the different behavior between the matrix-free version of the method of polarized traces, and its bi-directional variant. We can clearly observe that the execution time increases faster for the original method of polarized traces, which is mainly due to the increase of the number of iteration to convergence. This trend is more marked in the case of the positive perturbation as explained in Section \ref{section:ABC}. In particular, we observe that the execution time (due to an increase in the number of iterations) grows as $\cO(N \omega^{1/2})= \cO(N^{5/4} )$ instead of $\cO(N\log{N})$  and thus increasing overall the complexity of the global algorithm. In addition, we can observe that the constants are lower for the bi-directional preconditioner. We point out that the bottleneck of both preconditioners is the local solves, which are performed using a multifrontal solver; the main difference between both methods is the number of such solves at each iteration: for the method of polarized traces we need $7$ solves\footnote{The larger amount of solves needed for the matrix-free version of the method of polarized traces is due to the application of the integral operator using local solves at each GMRES iteration, in addition to the solves required for the application of the preconditioner.} per layer and per iteration, whereas the bi-directional preconditioner needs only 4; in addition, the bi-directional preconditioner converges in a lower number of iterations.

\begin{figure}[h]
\begin{center}
\includegraphics[trim= 0mm 0mm 00mm 0mm, angle =0,clip, width=6.3cm]{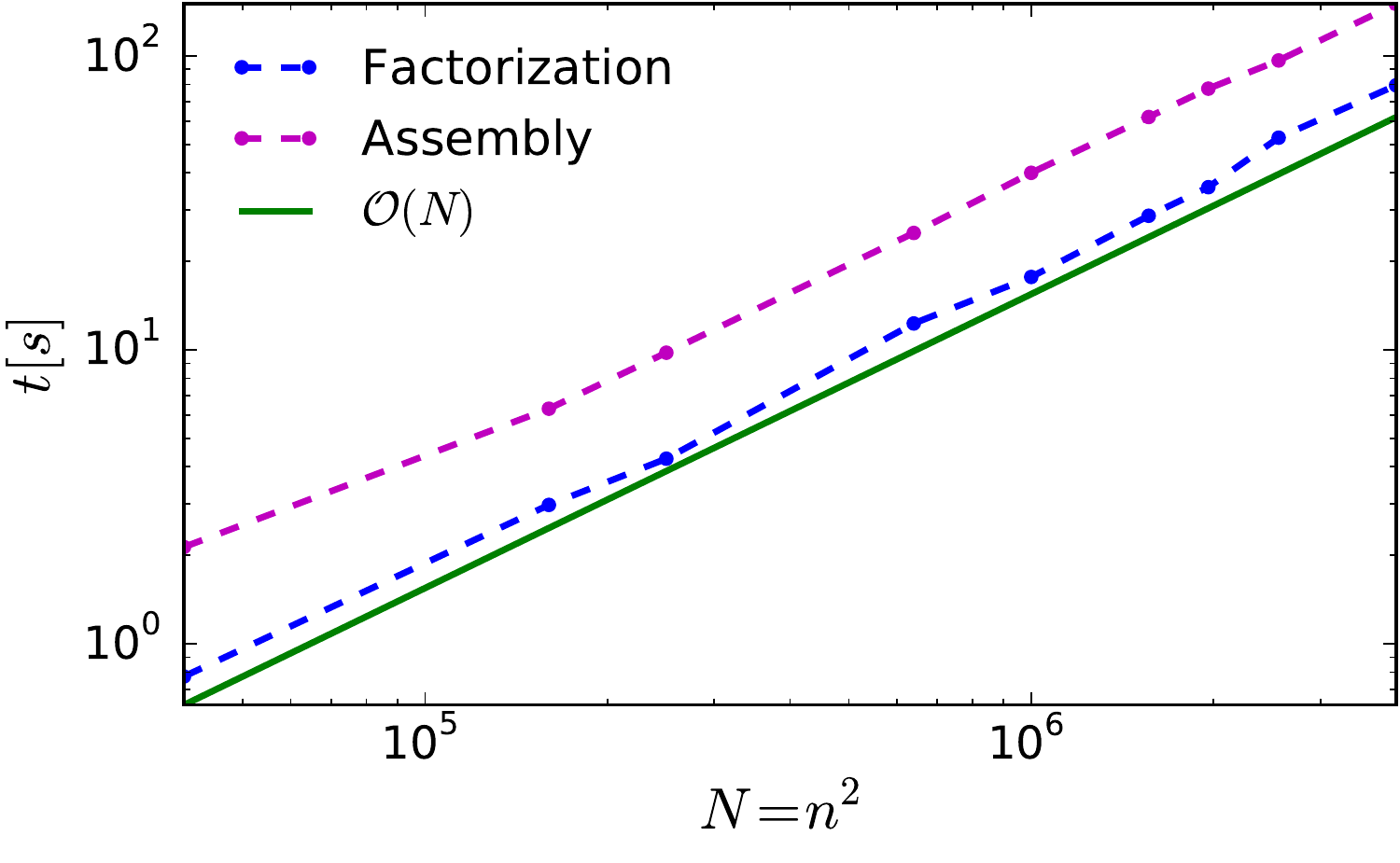}
\includegraphics[trim= 0mm 0mm 00mm 0mm, angle =0,clip, width=6.3cm]{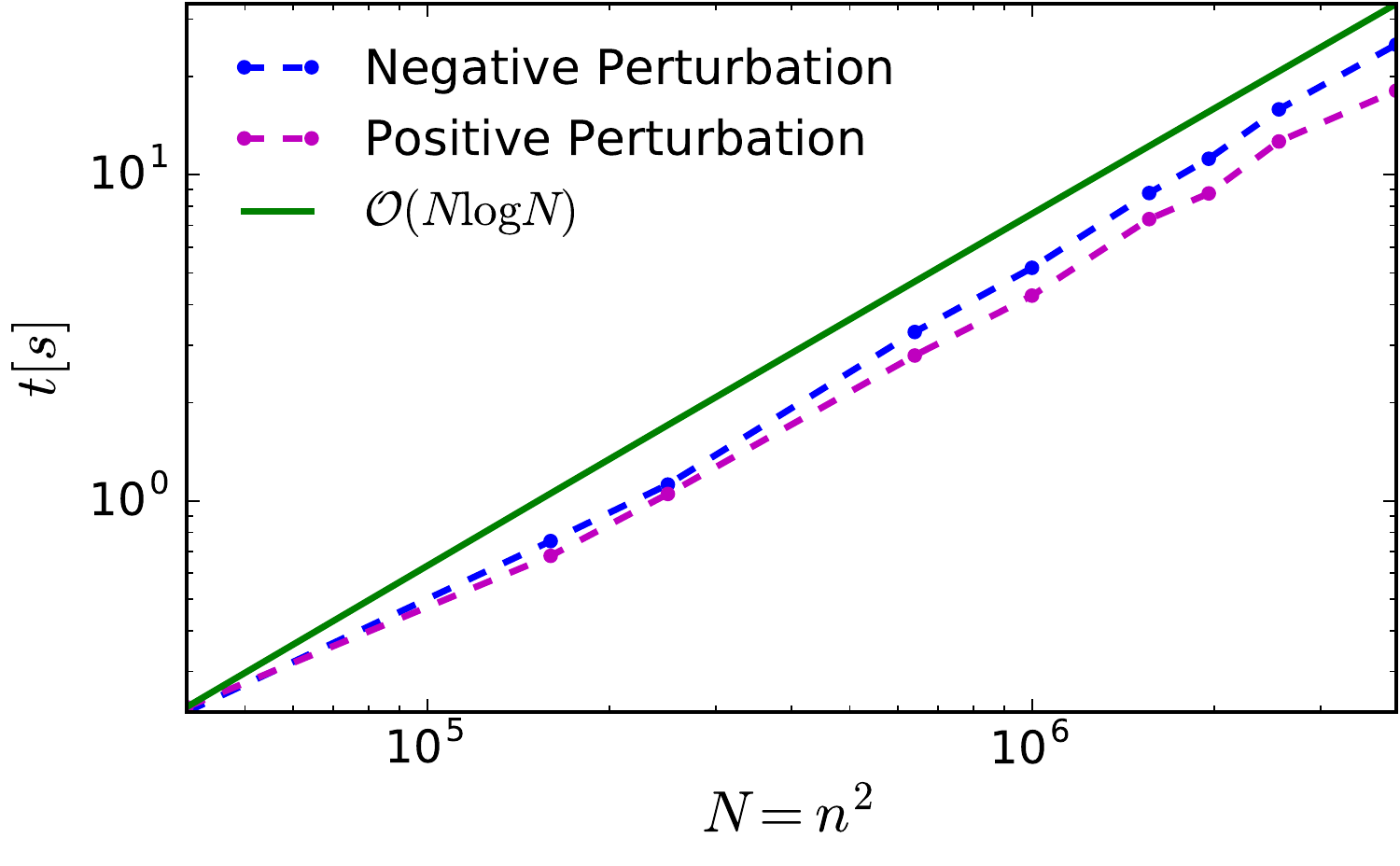}
\end{center}
\caption{{\it Left:} execution time of the off-line stage;  {\it Right:} execution time of the on-line stage of solving the Sparsified Lippmann-Schwinger equation using GMRES preconditioned  with Alg.~\ref{alg:birirectional_preconditioner}. The tolerance was fixed to $10^{-6}$, and  $\omega h = \cO(1)$.}\label{fig:scaling_sparsified_lippmann_schwinger}
\end{figure}

Fig.~\ref{fig:scaling_sparsified_lippmann_schwinger} ({\it left}), shows the scaling for the factorization of the local problems and the assembly of the sparsified system. In addition,  Fig.~\ref{fig:scaling_sparsified_lippmann_schwinger} ({\it right}), shows the execution time for the solution of the sparsified system using the bi-directional algorithm as reported in Table \ref{table:numerical_experiment_bi_directional} but in logarithmic scale. We can observe from Fig.~\ref{fig:scaling_sparsified_lippmann_schwinger} that the assembly and factorization cost is $\cO(N)$ and that the online cost of one inner iteration is $\cO( N \log{N})$, and, in particular, $n_{\mbox{iter}} = \cO(\log{\omega})$.

\subsection{Wave scattering}

\begin{figure}[h]
\begin{center}
\includegraphics[trim= 30mm 10mm 35mm 15mm, angle =0,clip, width=6.3cm]{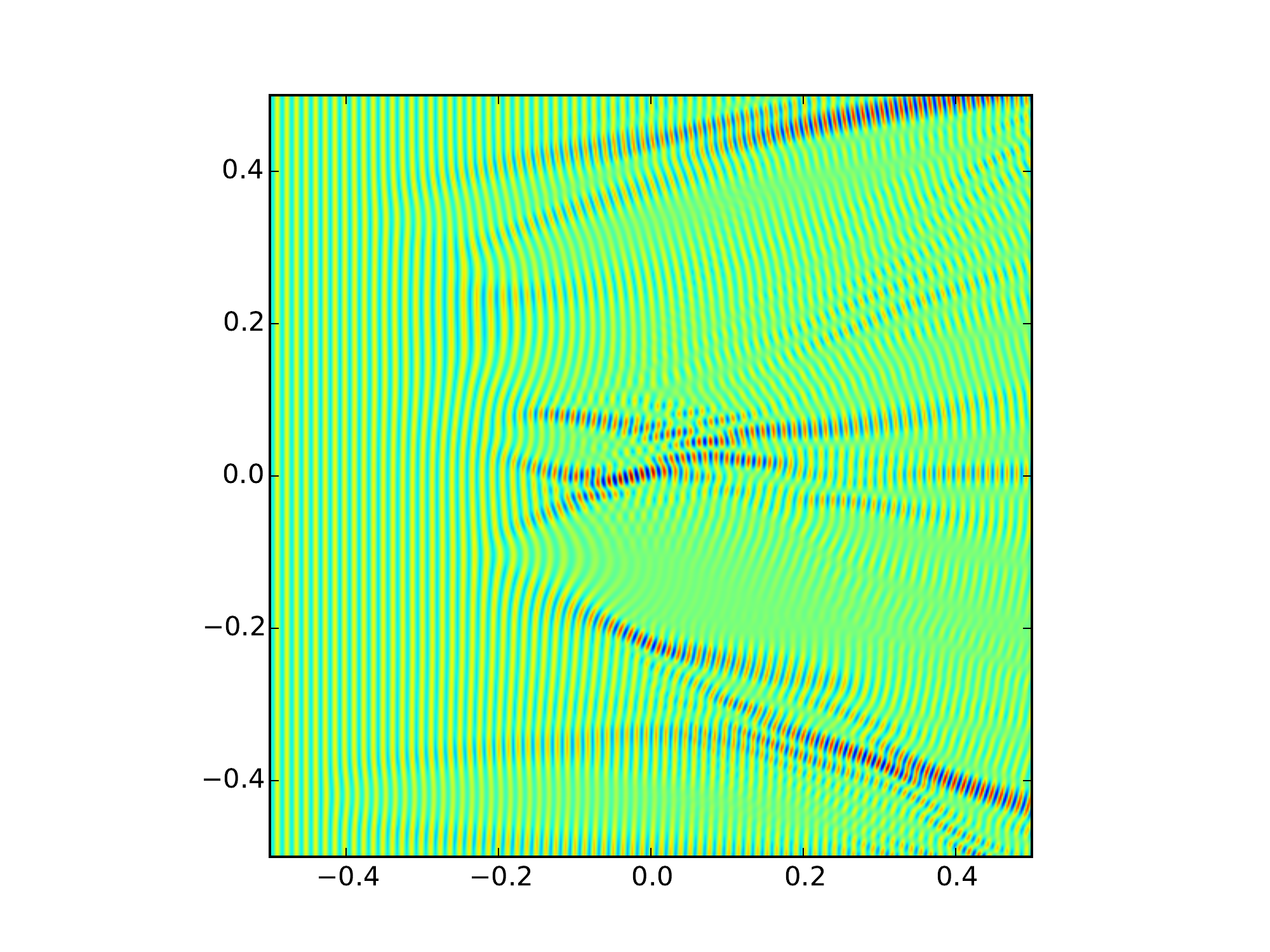}
\includegraphics[trim= 30mm 10mm 35mm 15mm, angle =0,clip, width=6.3cm]{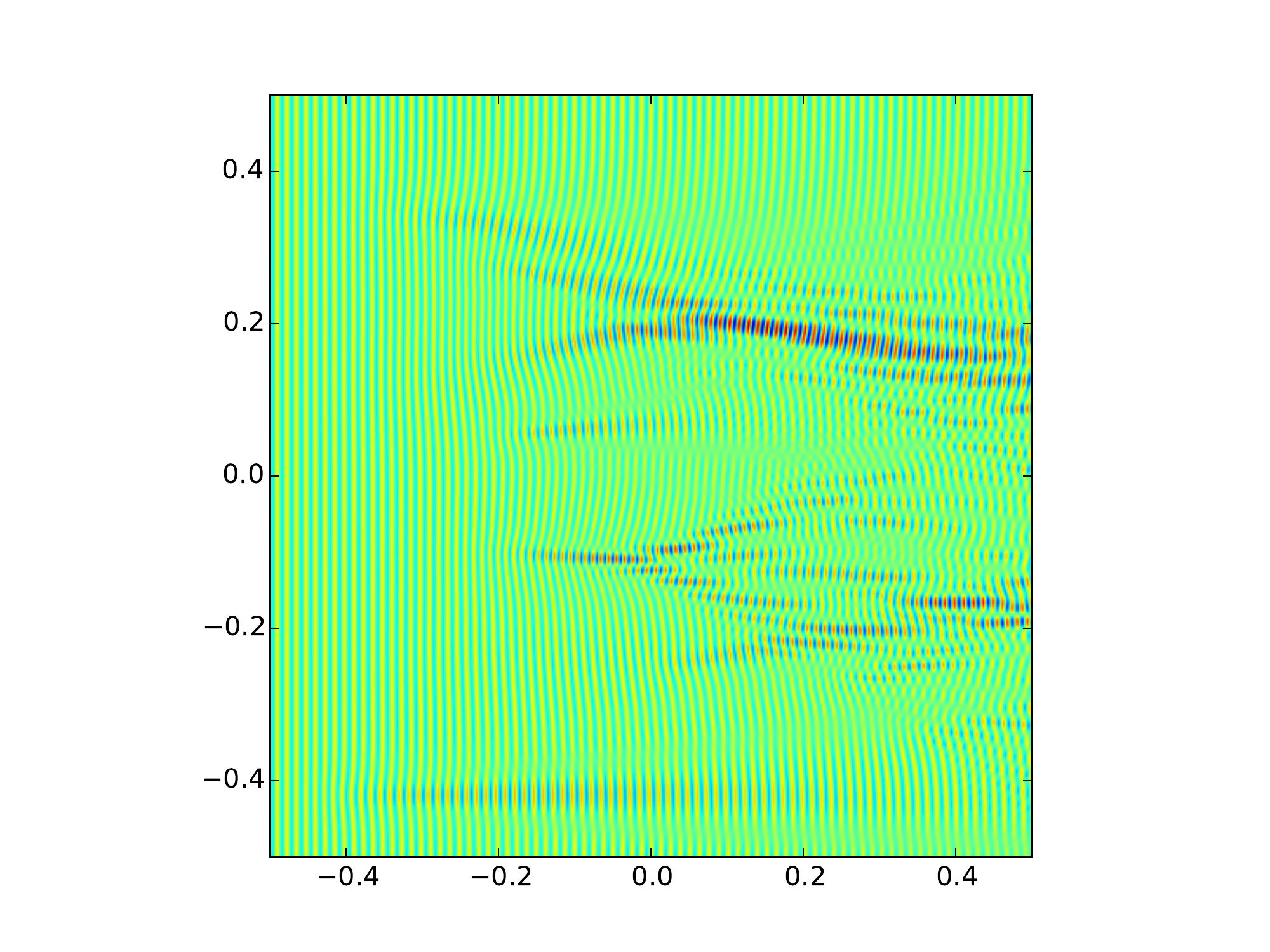}
\end{center}
\caption{Real part of the total wavefield generated by the scattering of a plane wave impinged on the negative perturbation ({\it left}) and on the positive perturbation ({\it right}), in both examples we have $\omega = 500$.}\label{fig:wavefields1}
\end{figure}

We test the performance of the two-level preconditioner by solving the Lippmann-Schwinger equation \eqref{eq:discrete_Lipmann-Schwinger} associated to three problems in wave scattering. We use three different perturbation profiles, the negative smooth perturbation, the positive smooth perturbation and Gaussian bumps as shown in Fig.~\ref{fig:scattering} ({\it left}) and ({\it right}). The profiles were chosen to depict the behavior of the method in three generic settings for wave propagation: a defocussing medium (negative perturbation), a focusing medium (positive perturbation) and a medium with multiple scattering (Gaussian bumps). We can observe the total wavefield scattered by a plane-wave impinging on each perturbation in Fig.~\ref{fig:wavefields1} ({\it left}), Fig.~\ref{fig:wavefields1} ({\it right}) and Fig.~\ref{fig:wavefields2} ({\it left}), respectively.

\begin{figure}[h]
\begin{center}
\includegraphics[trim= 0mm 0mm 00mm 0mm, angle =0,clip, width=6.3cm]{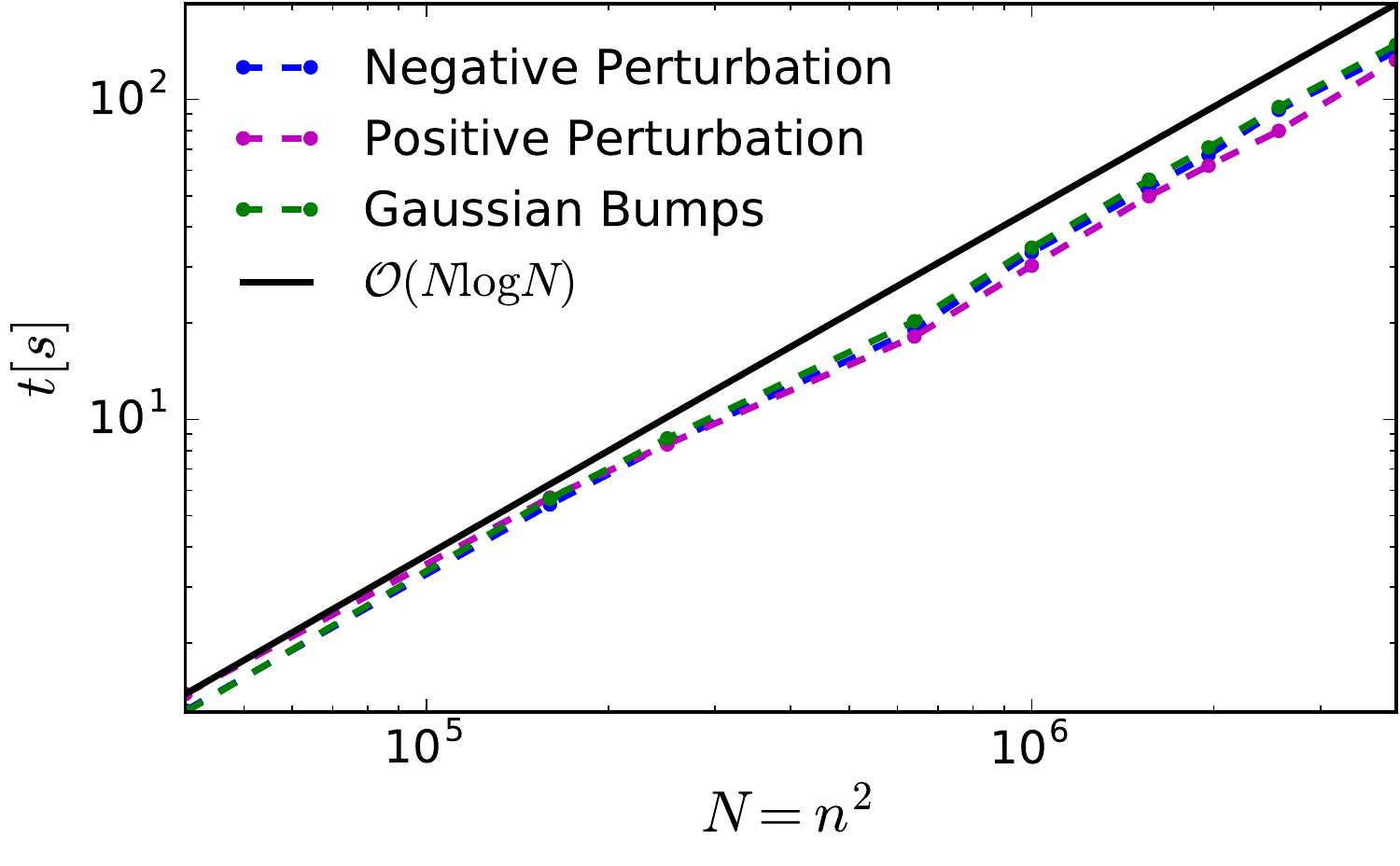}
\end{center}
\caption{ Average execution time for the on-line stage of solving the high-frequency Lippmann-Schwinger equation with GMRES (tolerance of $10^{-10}$) preconditioned with Alg.~\ref{alg:Lipmann-Schwinger_preconditioner} with an inner tolerance of $10^{-3}$, for the three different perturbation profiles; $\omega h = \cO(1)$.}\label{fig:scaling2}
\end{figure}

For each perturbation profile we solve a system of increasing size, keeping $\omega h$ constant. Each system is solved using GMRES with a tolerance of $10^{-10}$, preconditioned with Alg.~\ref{alg:Lipmann-Schwinger_preconditioner} with a inner tolerance of $10^{-3}$; the domain is decomposed in $L = n/50$ subdomains, and each subdomain is extended 10 grid points in each direction (see Fig.~\ref{fig:DDM_sketch}). At each frequency, we compute the wavefield generated by a set of 64 different incident waves impinging on the perturbation by solving the Lippmann-Schwinger equation, and we report the average execution time in Fig.~\ref{fig:scaling2}. From  Fig.~\ref{fig:scaling2} we can observe that the execution time scales as $\cO(N \log N)$ for each perturbation profile, which corroborates the claims made in the introduction. For the off-line stage we obtain the same execution times as in Fig.~\ref{fig:scaling_sparsified_lippmann_schwinger} ({\it left}).

\begin{figure}[h]
\begin{center}
\includegraphics[trim= 0mm 0mm 00mm 0mm, angle =0,clip, width=6.3cm]{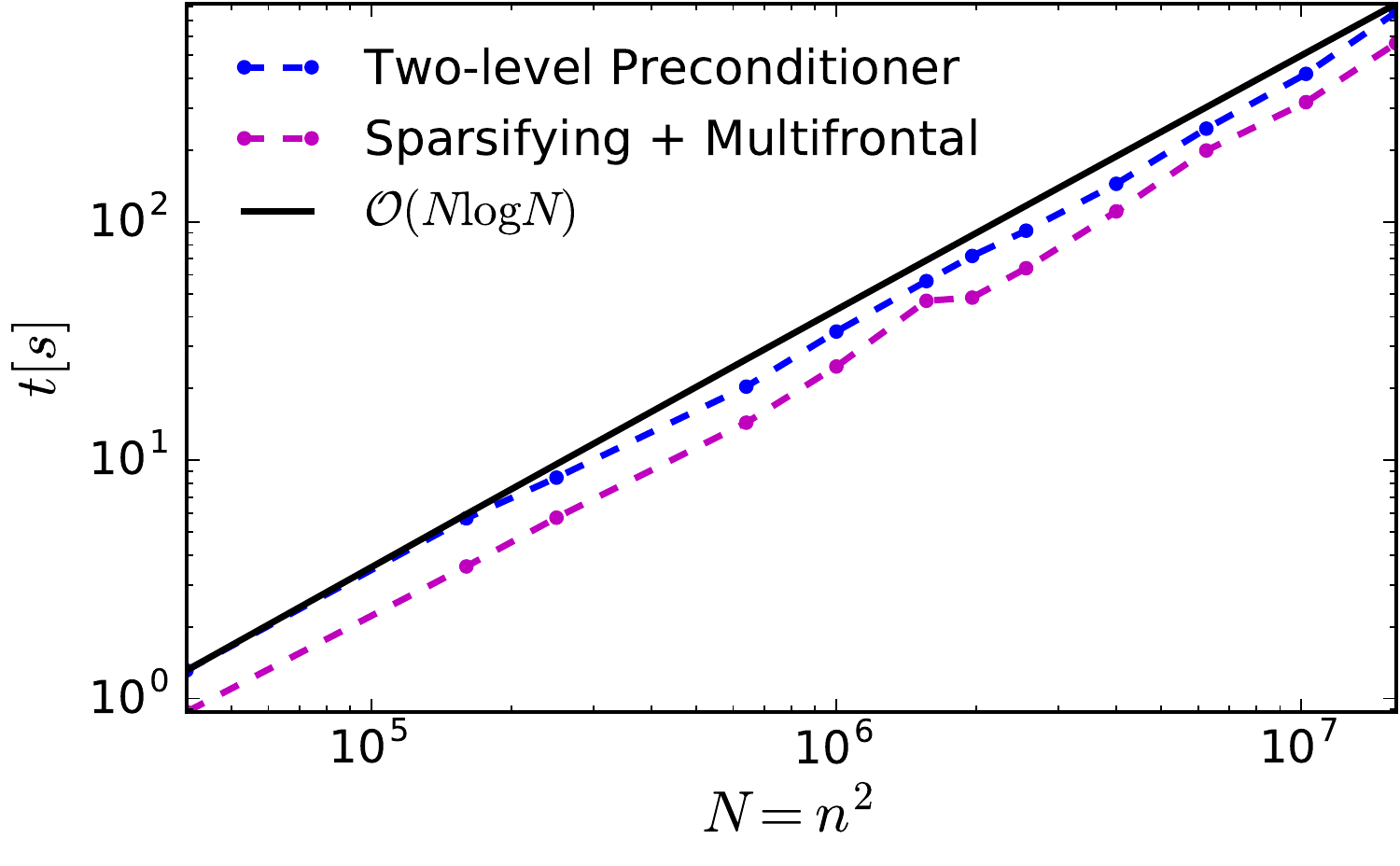}
\includegraphics[trim= 0mm 0mm 00mm 0mm, angle =0,clip, width=6.3cm]{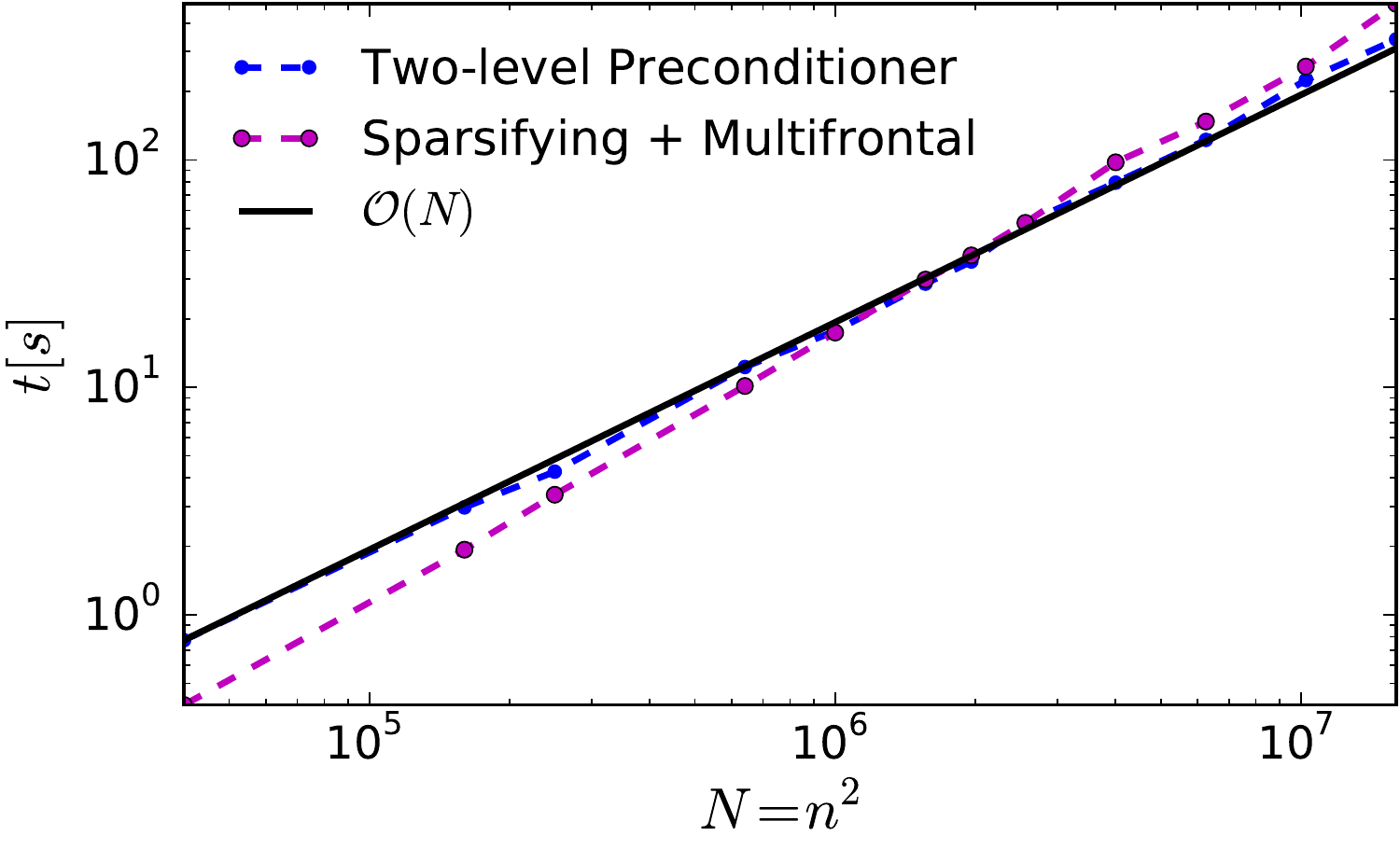}
\end{center}
\caption{ Left: average execution time for the on-line stage of solving the Lippmann-Schwinger with the Gaussian bumps perturbation with the original sparsifying preconditioner, i.e., using multifrontal methods, and using the two-level preconditioner with $\omega h = \cO(1)$; Right: execution time of the off-line stage for the two-level preconditioner and the original sparsifying preconditioner.}\label{fig:scaling3}
\end{figure}

Moreover, for the sake of comparison, we time the execution time of solving the Lippmann-Schwinger equation using the original sparsifying preconditioner \cite{Ying:Sparsifying_Preconditioner_for_the_Lippmann--Schwinger_Equation} (i.e., solving \eqref{eq:global_sparse_system} via a multifrontal solver \cite{Davis:UMFPACK}). For the comparison we used the multiple scattering perturbation in Fig.~\ref{fig:scattering} ({\it right}), with the same set of frequencies and incoming waves as before. We use the implementation of the sparsifying preconditioner in \cite{Zepeda:Lippmann_Schwinger}. The average execution times of the on-line stage are reported in Fig.~\ref{fig:scaling3} ({\it left}) together with the averages times for the two-level preconditioner. From Fig.~\ref{fig:scaling3} ({\it left}) we can observe that, unsurprisingly, the original sparsifying preconditioner has the same asymptotic on-line cost as the two-level preconditioner presented in this paper; however, the main difference lies in the off-line stage, in which the multifrontal solver has a complexity of $\cO(N^{3/2})$ compared to $\cO(N)$ for the two-level preconditioner as shown in Fig.~\ref{fig:scaling3} ({\it right}).

\subsection{Plasma physics}

\begin{figure}[h]
\begin{center}
\includegraphics[trim= 0mm 10mm 20mm 10mm, angle =0,clip, width=8cm]{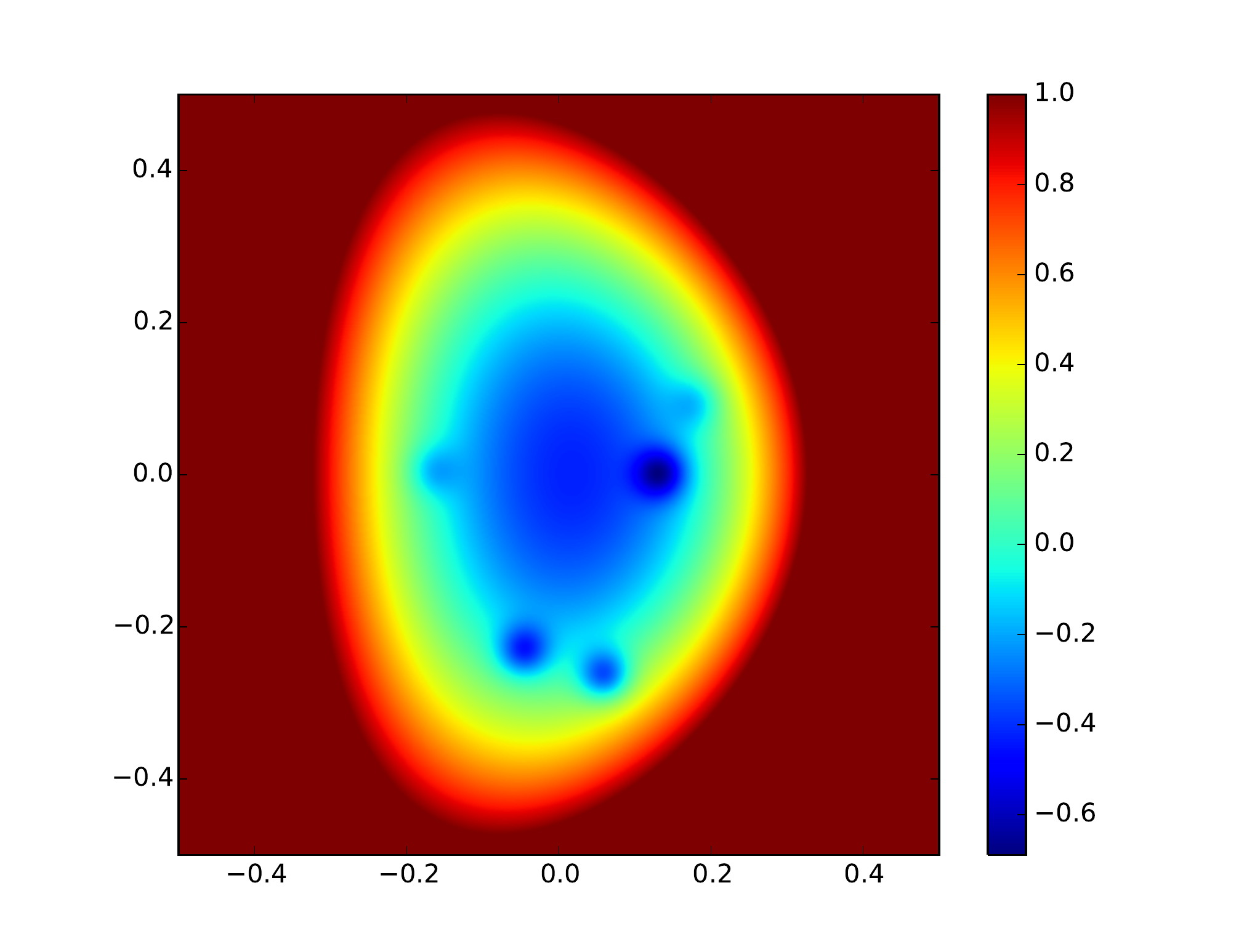}
\end{center}
\caption{Perturbation for the plasma physics example (taken from \cite{ambikasaran_greengard:Fast_adaptive_high_order_accurate_discretization_of_the_Lippmann-Schwinger_equation_in_two_dimension}) we can observe the at the interior the perturbation becomes negative, thus yielding a strongly elliptic problem locally, thus becoming an evanescent medium.}\label{fig:plasma}
\end{figure}

\begin{figure}[h]
\begin{center}
\includegraphics[trim= 30mm 10mm 35mm 15mm, angle =0,clip, width=6.3cm]{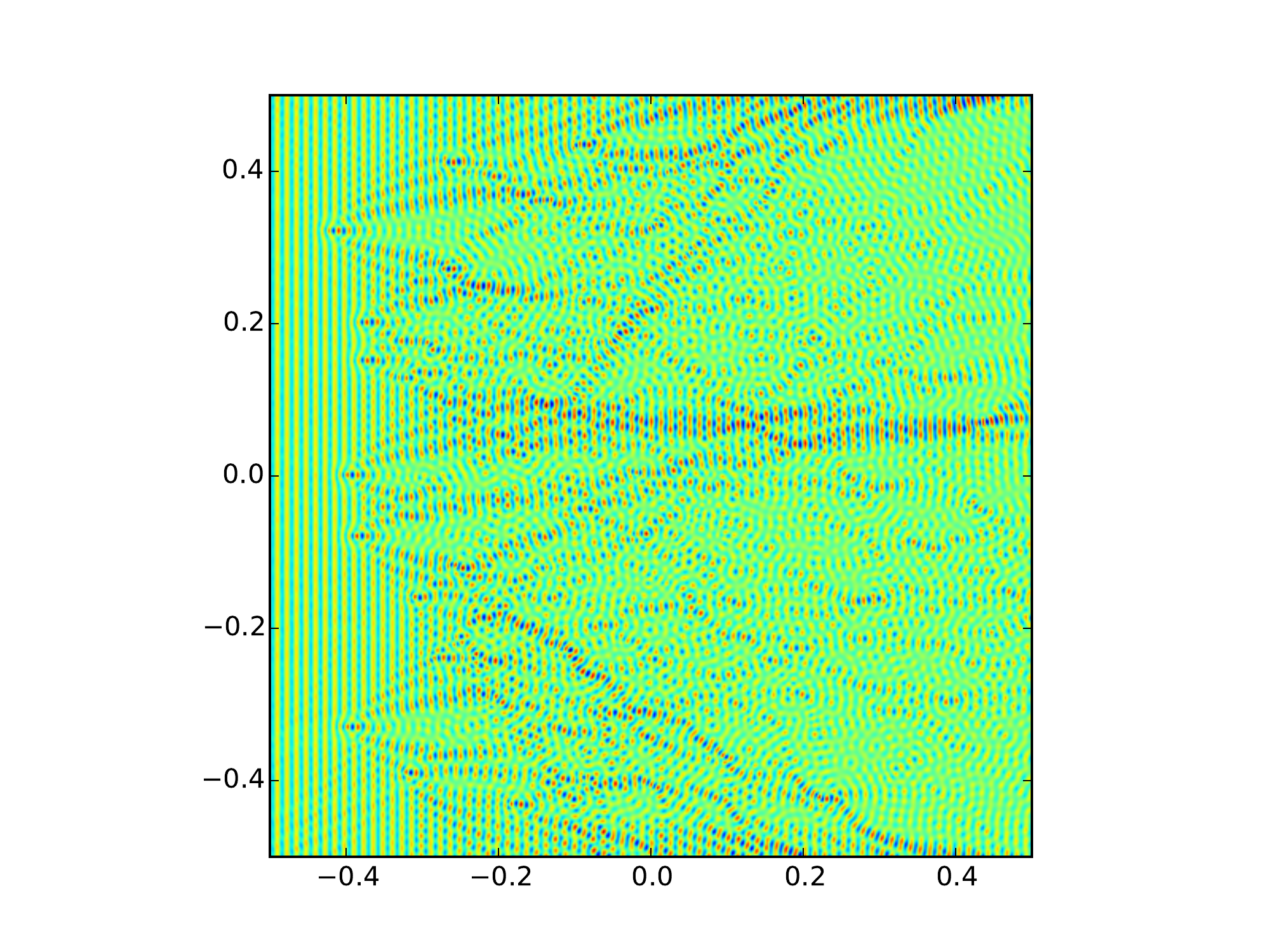}
\includegraphics[trim= 30mm 10mm 35mm 15mm, angle =0,clip, width=6.3cm]{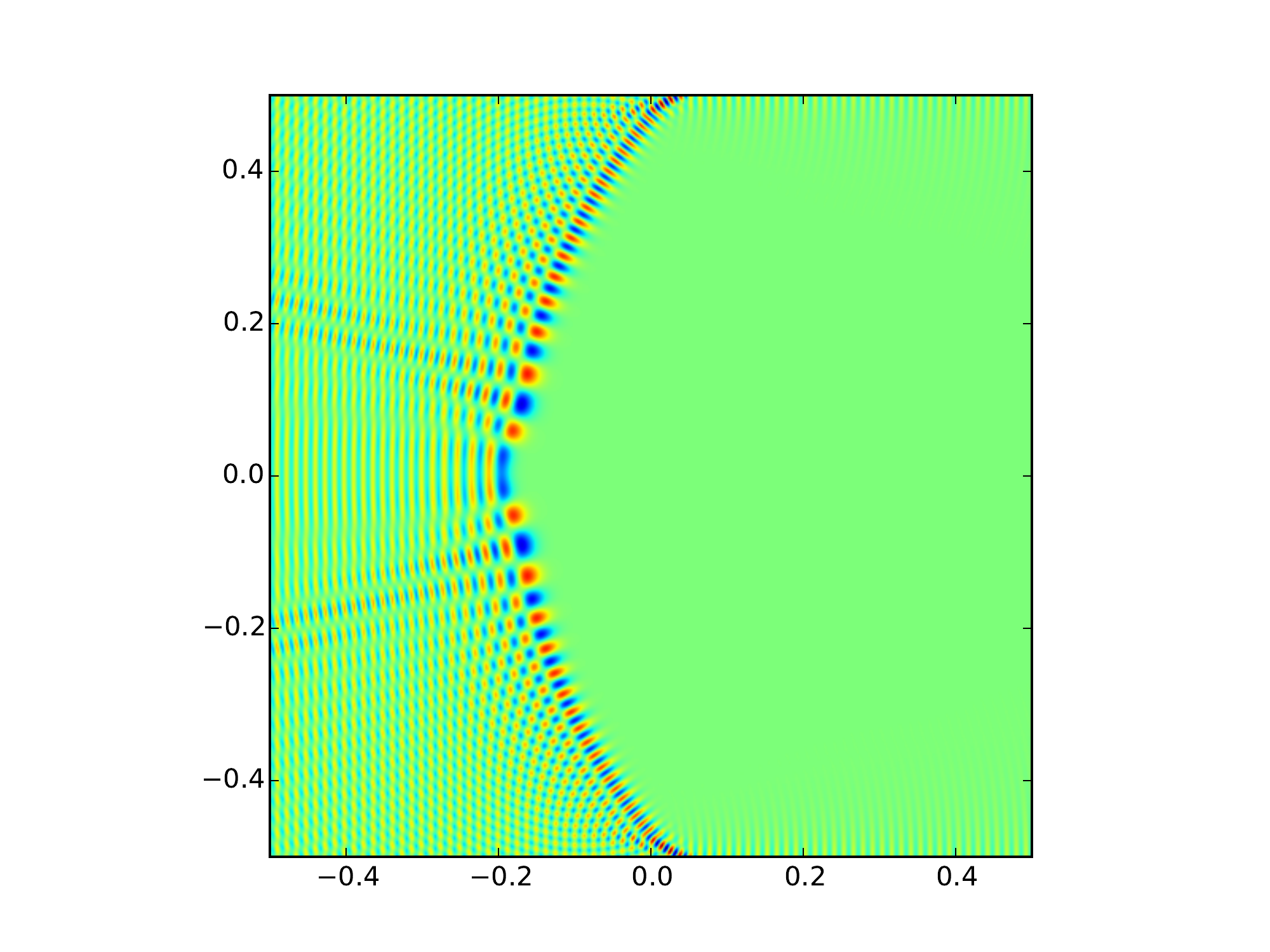}
\end{center}
\caption{Real part of the total wavefield generated by the scattering of a plane wave in the Gaussian bumps model (left) and in the plasma model (right). In both examples $\omega = 500$.}\label{fig:wavefields2}
\end{figure}

\begin{figure}[h]
\begin{center}
\includegraphics[trim= 0mm 0mm 00mm 0mm, angle =0,clip, width=6.3cm]{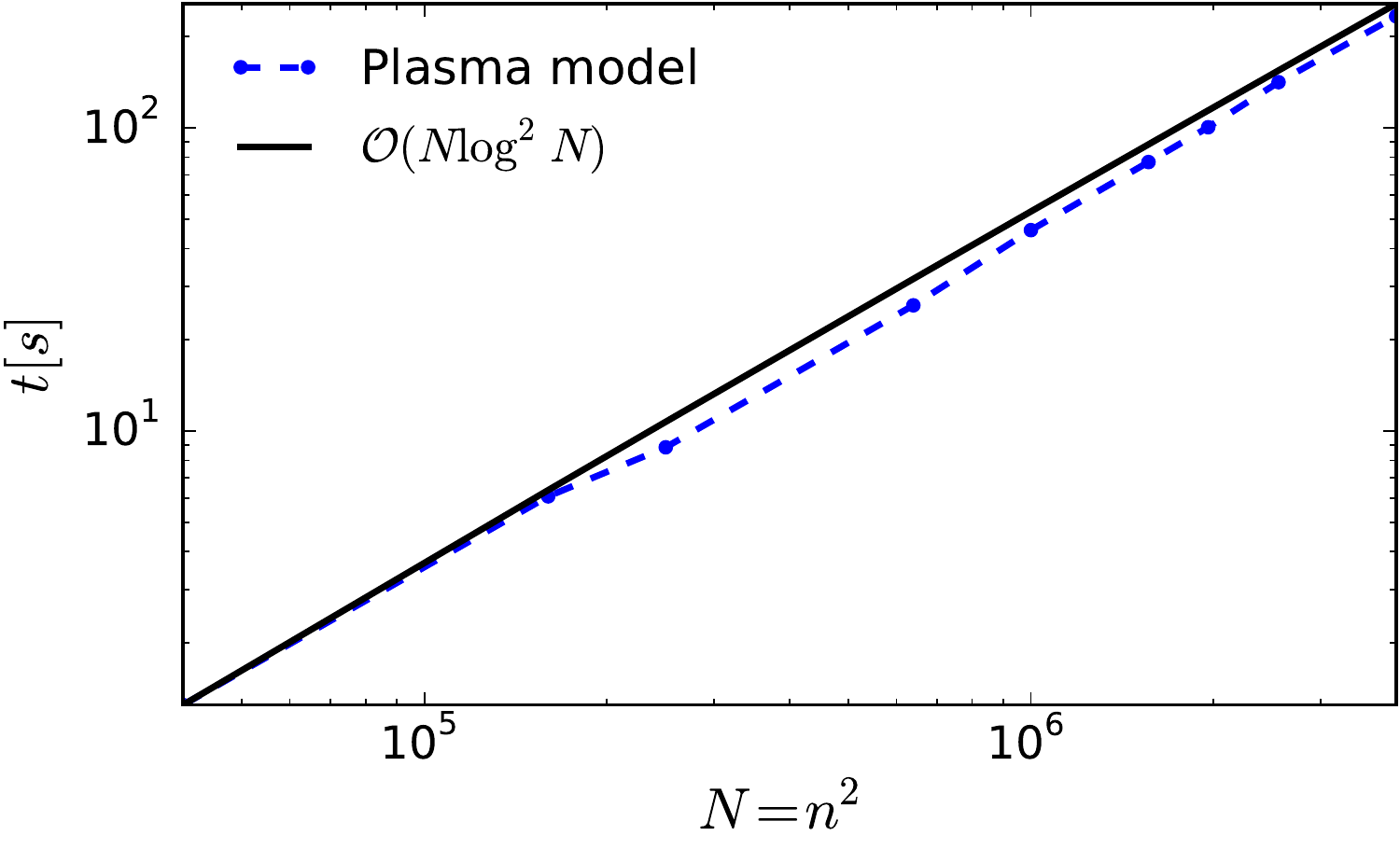}
\includegraphics[trim= 0mm 0mm 00mm 0mm, angle =0,clip, width=6.3cm]{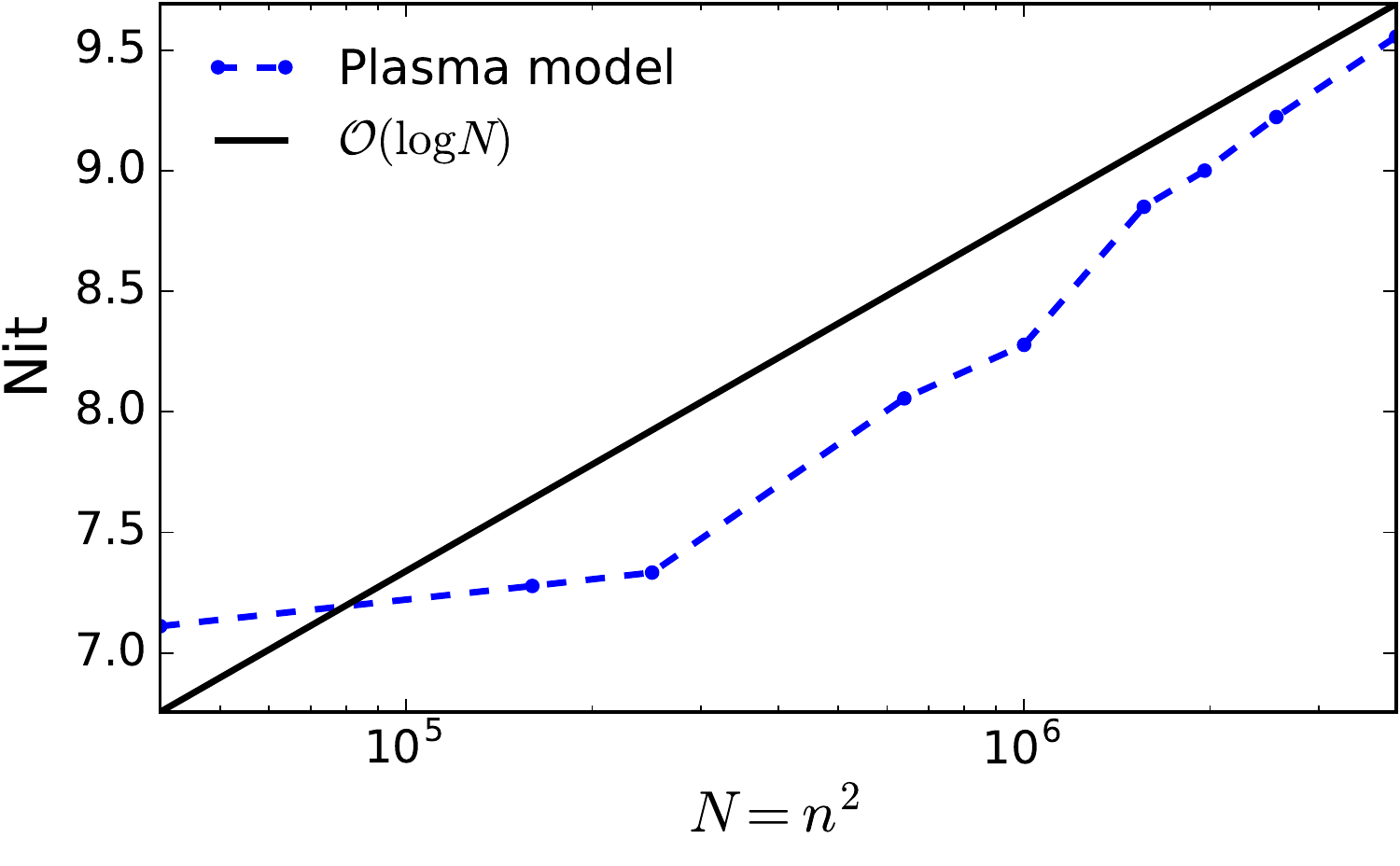}
\end{center}
\caption{Average execution time for solving the Lippmann-Schwinger equation using the plasma perturbation profile ({\it left}),  and the average number of outer iterations needed for convergence ({\it right}). The problem was solved with GMRES (tolerance of $10^{-10}$) preconditioned with the two-level preconditioner with an inner tolerance of $10^{-3}$ and  $\omega h = \cO(1)$. The extra logarithmic factor comes from the increasing number of outer iterations.}\label{fig:scalingplasma}
\end{figure}

Finally, we test the performance of the two-level preconditioner using a problem arising from plasma physics, in particular, a radio-frequency (RF) wave propagating in a fusion reactor. Under some simplifications, the propagation of a RF wave can be described using a cold plasma model \cite{Weitzner:Lower_hybrid_waves_in_the_cold_plasma_model}. We borrow from \cite{ambikasaran_greengard:Fast_adaptive_high_order_accurate_discretization_of_the_Lippmann-Schwinger_equation_in_two_dimension} a model\footnote{We refer the interested reader to Example 4.3 of \cite{ambikasaran_greengard:Fast_adaptive_high_order_accurate_discretization_of_the_Lippmann-Schwinger_equation_in_two_dimension} for further details.} corresponding to an idealized plasma density profile shown in Fig.~\ref{fig:plasma}. The physics in this case are different, an incoming wave propagates until it hits a cut-off region in which the problem becomes strictly elliptic or evanescent, thus causing the waves to reflect, as depicted in Fig.~\ref{fig:wavefields2} ({\it right}).

We solve a system of increasing size, keeping $\omega h$ constant; we decompose the domain such that we have roughly $5$ wavelengths inside each subdomain, and each subdomain is extended 10 grid points in each direction. At each frequency, we compute the wavefield generated by a set of 64 different incident waves impinging on the perturbation by solving the Lippmann-Schwinger equation and we report the average execution times in Fig.~\ref{fig:scalingplasma} ({\it left}).

In the case of the plasma example, we observe from Fig.~\ref{fig:scalingplasma} ({\it left}) a slightly bigger scaling. The main reason for this scaling is the fact that the number of outer iteration grows logarithmically, as shown in Fig.~\ref{fig:scalingplasma} ({\it right}), which is presumably caused by the transition of the physics at the cut-off region. Nonetheless, the number of inner iterations grows like $\cO(\log{N})$, the same as in the wave scattering case, resulting in an overall complexity of $\cO(N \log^2{N})$, which, to the our knowledge, compares favorably to most of current solvers.

\section{Conclusions}
We have presented a new fast iterative method for the high-frequency Lippmann-Schwinger equation in 2D. The method has an asymptotic cost of $\cO(N \log{N})$  for the wave scattering problems presented in this paper.
%Moreover, for waves impinging in a cold plasma the complexity slightly deteriorates to $\cO(N \log^2{N})$.

Some future work and extensions are: the design of more efficient ABC, in order to use only one directional preconditioner, thus decreasing the memory footprint; extensions to 3D, in which a nested approach can be used to obtain low complexity algorithms; parallelization of the solvers using state-of-the-art linear solvers; reduction of the constant via compressed linear algebra, which can provide extra savings in the plasma example, where the local ellipticity of the problem can be exploited.

Another line of research is to further reduce the on-line complexity by using the original reduction to a SIE of the method of polarized traces, coupled to the precomputation of integral operators and their compression.

We point out that it is envisageable to change the outer GMRES loop with the flexible GMRES \cite{Saad:A_Flexible_Inner-Outer_Preconditioned_GMRES_Algorithm} algorithm, and reduce the inner-level tolerance to accelerate the algorithm, without a penalty on the asymptotic complexity.

\section{Acknowledgments}

We would like to thank Lexing Ying, and Leslie Greengard for fruitful discussions; Antoine Cerfon for his help with the plasma example, and Laurent Demanet for computational resources.

% For the wave propagation we used the perturbation as in, with an incoming wave-field in the y direction

% \subsection{Waves in Plasma}

% We used the same example in

% \section{Conclusion}

%\bibliographystyle{siamplain}
\bibliography{GRF_integral_formulations.bib}

\end{document}